\documentclass[12pt,leqno]{article}
\usepackage{amssymb,amsmath,theorem,latexsym,epsfig}
\usepackage{color}
\usepackage[all]{xy}
\usepackage{psfrag}
\usepackage{amsfonts}
\usepackage{amssymb}
\usepackage{mathrsfs}
\usepackage{graphicx}

\theoremstyle{plain}

\newtheorem{theo}{Theorem}[section]

\newtheorem{cor}[theo]{Corollary}

\newtheorem{conn}{Conjecture}[section]
\newtheorem{defi}{Definition}[section]
\newtheorem{prop}[theo]{Proposition}

\newtheorem{claim}[theo]{Claim}
\newtheorem{prop-defi}[theo]{Proposition-Definition}
\newtheorem{lemma-defi}[theo]{Lemma-Definition}
{\theorembodyfont{\rmfamily} \newtheorem{rem}{Remark}[section]}               
{\theorembodyfont{\rmfamily} \newtheorem{ex}{Example}[section]}
{\theorembodyfont{\rmfamily} }

\newlength{\picwidth}
\newlength{\miniwidth}
\newlength{\nanowidth}
\newlength{\melowidth}
\newlength{\maxwidth}
\DeclareFontFamily{U}{rsf}{}
\DeclareFontShape{U}{rsf}{m}{n}{
  <5> <6> rsfs5 <7> <8> <9> rsfs7 <10->  rsfs10}{}
\DeclareMathAlphabet{\mathscr}{U}{rsf}{m}{n}

\newcommand{\op}[1]{\operatorname{#1}}

\newcommand{\C}{{\mathbb C}}

\newcommand{\coh}{\op{Coh}}

\newcommand{\bC}{\mathbb{C}}

\newcommand{\ZZ}{{\mathbb Z}}
\newcommand{\CP}{{\mathbb{P}}}
\newcommand{\bw}{{\sf w}}
\newcommand{\bfw}{\bw}
\newcommand{\bfS}{{\sf S}}
\newcommand{\bfY}{{\sf Y}}
\newcommand{\bfq}{{\sf q}}
\newcommand{\mir}{\boldsymbol{|}}

\def\shE{{\cal E}}

\newcommand{\E}{\mathcal{E}}
\newcommand{\db}{D^{b}}

\newcommand{\supp}{\operatorname{supp}}

\newcommand{\perf}{\operatorname{{\sf perf}}}
\newcommand{\Hom}{\operatorname{Hom}}
\newcommand{\dsing}{D^{b}_{\op{sing}}}
\newcommand{\bP}{\mathbb{P}}
\newcommand{\bZ}{\mathbb{Z}}
\newcommand{\rhom}{R\Hom}
\newcommand{\ext}{\op{Ext}}
\newcommand{\Fuk}{{D\sf{Fuk}}}
\newcommand{\FS}{{D\sf{FS}}}
\newcommand{\bfa}{\boldsymbol{a}}
\newcommand{\bfb}{\boldsymbol{b}}
\newcommand{\bff}{\boldsymbol{f}}

\newcommand{\tn}[1]{\textnormal{#1}}

\begin{document}
\noindent

\title{Homological Mirror Symmetry for manifolds of general type}  
\author{A. Kapustin, L.Katzarkov \footnote{Partially supported by NSF
    Grant DMS0600800 and by Clay Math Institute}\,\,,  D. Orlov,
  M. Yotov} 
\date{ }
\maketitle

\begin{flushright}
\begin{minipage}[c]{2.5in} {To Fedya, our  teacher and friend, with
    admiration} 
\end{minipage}
\end{flushright}

\bigskip

\tableofcontents

\section{Introduction}

This paper discusses Homological Mirror Symmetry for manifolds of
general type and the way this symmetry interacts with interesting
conjectural dualities in complex and symplectic geometry.
Traditionally, \cite{candelasetal,cox-katz}, the mirror symmetry
phenomenon is studied in the most relevant for the Physics case of
Calabi-Yau manifolds. It was observed that some pairs of such
manifolds ("mirror partners"), which naturally appear in mathematical
models for physics string theory, exhibit properties (have numerical
invariants) that are symmetic ("mirror") to each other. There were
several attempts to mathematically formalize this mirror phenomenon,
one of which was the celebrated Homological Mirror Symmetry Conjecture
stated by Kontsevich in 1994. Given two Calabi-Yau mirror partners,
they are naturally complex manifolds and symplectic
manifolds. Accordingly, there are two categories defined on each of
them: the bounded derived category of coherent sheaves (reflecting the
complex structure of the manifolds and "physically" related to
D-branes of type B), and the derived Fukaya category (a symplectic
construct and related to D-branes of type A). What the Conjecture
states is that, for any pair of mirror partners, the category of
$B$-branes on each of the manifold is equivalent to the category of
$A$-branes on the mirror partner.

Kontsevich's conjecture was verified for large number of Calabi-Yau
manifolds: elliptic curves (Polischuk-Zaslow) \cite{PoZa}, abelian varieties
(Fukaya) \cite{fukaya}, K3 surfaces (Seidel) \cite{SAI1}.

A basic ingredient in these considerations is the construction of
(physics) mirror partners. An approach by Hori and Vafa, interpreting
the latter as a holomorphic fibration (Landau-Ginzburg model), turned
out to be working well even beyond the case of Calabi-Yau
manifolds. Specifically, this approach provides "mirror partners" for
all complete intersections in toric manifolds. When the pairs of
manifolds constructed this way are not Calabi-Yau, there is no
physical ground to call them mirror symmetric. Moreover, the
Landau-Ginzburg partner in this case, as given by the Hori-Vafa
approach, is not compact. So, it was expected that appropriate
modifications had to be made in the construction of this partner as
well as in the definition of the relevant categories in order to make
it possible to extract valuable geometric information out of this
fibration. Furthermore, the formulation of the conjecture itself
needed to be changed in order to make it a true statement. Suppose the
$A$- and $B$-brane categories can be defined (with modifications on the Landau-Ginzburg
partner) and compared. Then by definition, when these are
cross-equivalent as in the HMSC, the manifolds are called
(homologically) mirror symmetric. The task now would be to show that
any complete intersection in a toric manifold and its Landau-Ginzburg
partner provide a pair of (homologically) mirror partners. Indeed, an
of idea how to change the categories exists, and it was proved that
when the starting manifold is a del Pezzo surface (a Fano manifold),
then the Landau-Ginzburg partner is a (homologically) mirror partner of the
surface. The same result holds true for weighted projective planes and
Hirzebruch surfaces as well.

In this paper we discuss how the HMSC should be modified so that it
holds true in the case of manifolds of general type (with negative
first Chern class). The conjectures we make are based on numerous
examples some of which we present here. In particular, we study in
considerable detail the Landau-Ginzburg partner of the hyperelliptic
curves or of the associated (see \cite[Theorem~2.7]{BO}) intersections
of two quadrics.

Our main objective is to study how the HMSC for non-Calabi-Yau
manifolds will interact with interesting geometric dualities on the
$A$ and the $B$ side of the mirror correspondence. In this direction
we propose a conjecture relating the complex or symplectic geometry of
a pair of manifolds, one of which is a Fano, and the other is of
general type:

\

\medskip

\noindent
{\bfseries Conjecture\ref{con:side.A.equiv}} \ {\em Let $F$ and $G$ be
  projective manifolds. Suppose $F$ is a Fano manifold, and $G$ is a
  manifold of general type or a Calabi-Yau manifold. Suppose further
  that there exists a fully faithful functor 
\[
\Phi : D^{b}(G) \hookrightarrow
  D^{b}(F)
\] 
between the derived categories of $F$ and $G$. Write
  $\underline{F}$ and $\underline{G}$ for the $C^{\infty}$ manifolds
  underlyng $F$ and $G$. Then:
\begin{itemize}
\item[(i)] For every complexified K\"{a}hler class $\alpha_{F}$ on
  $F$ there exists a complexified  K\"{a}hler class $\alpha_{G}$ on
  $G$ and a fully faithful functor 
\[
\Psi : \overline{\Fuk}(\underline{G},\alpha_{G}) \hookrightarrow
\overline{\Fuk}(\underline{F},\alpha_{F})
\]
between the corresponding Karoubi completed Fukaya categories.
\item[(ii)] If $K_{\Phi} \in D^{b}(G\times F)$ is a kernel object for $\Phi$
  (such an object exists by \cite{orlov-K3}), then $\Psi$ is given by a kernel
  object $K_{\Psi}$in $\overline{\Fuk}(\underline{G}\times
  \underline{F},-p_{G}^{*}\alpha_{G} + p_{F}^{*}\alpha_{F})$ which is
  uniquely determined by $K$.
\end{itemize}
}

\

\bigskip

\noindent
In particular, this conjecture predicts the existence of a functor (of
geometric nature) between the Fukaya categories of a curve of genus
two and the intersection of two quadrics in the projective
5-space. Using Seidel's recent proof \cite{SEID} of HMS for the genus
two curve, we can recast the above conjecture as a statement about
$B$-branes in the mirror Landau-Ginzburg models. Our main result,
Theorem~\ref{theo:main} is that in a large complex structure limit a
fully faithful functor between these Landau-Ginzburg models does
exist. This proves the large volume case of
Conjecture~\ref{con:side.A.equiv} when $G$ is a curve of genus two and
$F$ is the associated degree four Fano variety in $\mathbb{P}^{5}$.

The analogue of Theorem~\ref{theo:main} is expected to be true for $G$
a hyper-elliptic curve of genus $g$ and $F$ the associated
intersections of two quadrics in $\mathbb{P}^{2g+1}$. The recent proof
\cite{efimov} of mirror symmetry for higher genus curves by Efimov
again allows us to reformulate the conjecture as a statement about
holomorphic Landau-Ginzburg models and we expect that our method of
proof will work in that case as well.

The paper is organized as follows. In section~\ref{sec:definitions} we
collect all the changes needed in the classical (regarding Calabi-Yau
manifolds) definitions for the case of Landau-Ginzburg partners
(fibrations over the complex line). In particular, we discuss the
issue of the non-uniqueness of (partial) compactification and further
desingularization of the originally given Landau-Ginzburg partner. We
conclude the section by giving a warm-up example (of a degree four del
Pezzo surface) where our procedures work perfectly in agreement with
the expectations. In section \ref{sec-main}, based on physics
considerations, we state the main conjectures about how the HMS
conjecture should look like in the case of manifolds of general
type. In section \ref{sec-hyperelliptic} we study the Landau-Ginzburg
partners of hyper-elliptic curves. The emphasis here is on the case
of genus two curves.  In the following section \ref{sec-quadrics}, we
investigate the Landau-Ginzburg partner of the complete intersection
of two quadrics in the 5-dmiensional projective space, and compare it
with the results for genus 2 curve in section
\ref{sec-hyperelliptic}. In section \ref{sec:large} we compare the two
mirrors and prove a large volume version of
Conjecture~\ref{con:side.A.equiv}  modulo the Homological Mirror
Symmetry Conjecture for the intersection of two quadrics. In
section \ref{sec:flops} we show that the derived categories of three
dimensional Landau-Ginzburg models are preserved by simple flops, a
fact that is used repeatedly in our constructions of good models
for Landau-Ginzburg mirrors. In the last section \ref{sec:match} we
try to make the mirror correspondence for genus two curves more
explicit by matching explicitly representattive objects in the Fukaya
category with objects in the mirror Landau-Ginzburg model and
comparing the corresponding Floer homologies and Ext
groups.

Before we proceed, let us state the conjectures explaining how the HMS
should behave in the case of manifolds of general type. To explain the
conjectures we will need to introduce some notations. Consider
\begin{description} 
\item[\fbox{$(\underline{X},\omega)$} \ :] a symplectic manifold
associated with a projective manifold $X$ of general type,
i.e. $\underline{X}$ is the $C^{\infty}$-manifold underlying $X$ and
$\omega$ is a symplectic form such that $c_{1}(X) = [-\omega]$. For
instance take $\omega$ to be the curvature of a Hermitian connection
on the canonical line bundle $K_{X}$.
\item[\fbox{$\bw : Y \to \mathbb{C}$} \ :] a holomorphic
Landau-Ginzburg mirror of the symplectic manifold
$(\underline{X},\omega)$. For instance, in the case when $X$ is a
complete intersection in a toric variety, the mirror $\bw : Y \to
\mathbb{C}$ can be constructed by the Hori-Vafa prescription
\cite{HV}. 
\end{description}

\

\noindent
With this notation we have the following:

\

\noindent
{\bfseries Conjecture~\ref{con-AtoB}:} {\em
 If $(\underline{X},\omega)|(Y,\bw)$ is such a mirror
  pair, then:
\begin{itemize} 
\item The derived Fukaya category
$\Fuk(\underline{X},\omega)$ embeds as a direct summand into the
category $D^b(Y,\bw)$ of $B$-branes of $\bw : Y \to
\mathbb{C}$.
\item The orthogonal complement of
$\Fuk(\underline{X},\omega)$ in $D^b(Y,\bw)$ is very simple: it is a
direct sum of several copies of the category of graded modules over a
Clifford algebra of a symmetric bilinear form on a complex vector
space of dimension $n=\dim_{\mathbb{C}} Y$.
\end{itemize}
}

\

\medskip

\noindent
There is also an analogous mirror conjecture in which the $A$ and $B$
sides of the theory are switched. More precisely if denote the
$C^{\infty}$-manifold underlying $Y$ by $\underline{Y}$, then the
complex manifold $X$ of general type should determine a symplectic
structure $\eta$ on $\underline{Y}$ with respect to which $\bw :
\underline{Y} \to \mathbb{C}$ becomes a symplectic Lefschetz fibration.

\

\noindent
With this notation we have the following:

\

\medskip

\noindent
{\bfseries Conjecture~\ref{con-BtoA}:} {\em
If $X|(\underline{Y},\bw,\eta)$ is such a mirror
  pair, then the category $D^{b}(X)$ of $B$-branes on $X$ is
  equivalent to the $A$-brane category $\FS(\underline{Y}_{D},\bw,\eta)$
  of a potential $\bw : \underline{Y}_{D} \to D$, where $0 \in D \subset
  \mathbb{C}$ is a suitably chosen disk, and $\underline{Y}_{D} =
  \bw^{-1}(D)$. 
}                           

\

\bigskip

\noindent
{\bfseries Remark:} \ While this work has been written two related
works have appeared. The first is Seidel's paper \cite{SEID} in which
he proves HMS of genus 2 curves by a different approach.  The second
is the recent work of A.Efimov \cite{efimov} who had proved the HMS
conjecture for curves of arbitrary genus.

\

\medskip

\noindent
{\bf Acknowledgments:} We have benefited greatly from discussions
M.Kontsevich, P.Seidel, K.Hori, C.Vafa and E.Witten.  We are very
thankful to D.Auroux and T.Pantev, without whom this paper would have
never been written.  We would like to thank IPAM, UCLA and especially
M.Green and H.D.Cao for the wonderful working conditions of the
program Symplectic Geometry and Physics, during which the work on this
paper was initiated. The first two authors are grateful to KITP, Santa
Barbara and especially D.Gross, R.Dijgraaf, D.Freed and T.Pantev, for
the wonderful athmosphere during the program Geometry, Duality and
Strings, when most of this work was done.

\section{Some Definitions} \label{sec:definitions}

We will use the setup of HMS in situations when one of the partners is
a Landau-Ginzburg fibration.  For projective manifolds the categories
of topological branes are well understood \cite{mirrorbook,KO}. In the
simplest setup the category of branes in the $B$ model is identified
with the bounded derived category of coherent sheaves $D^{b}(X)$ on a
complex manifold $X$, and the category of $A$-branes is the derived
Fukaya category $\Fuk(\underline{Y},\omega)$ of a compact symplectic
manifold $(\underline{Y},\omega)$ \cite{OOO}. It is a triangulated
category which is the homotopy category of the category of twisted
complexes over the geometrically defined Fukaya category of Lagrangian
submanifolds equipped with complex local systems. More generally one
has to enhance the categories $D^{b}(X)$ and
$\Fuk(\underline{Y},\omega)$ with natural dg or
$A_{\infty}$-structures.

In the Landau-Ginzburg setting the categories of braness have to be
modified in an appropriate manner \cite{mirrorbook,KO}.
We begin by recalling the definition of
the categories of D-branes in topological Landau-Ginzburg theories.
First we deal with the symplectic ($A$-brane) side of the picture.  We
will follow Seidel's treatment of the $A$-twist of the Landau-Ginzburg
model \cite{seidel-book}. Historically the idea was introduced first
by M.Kontsevich and later by K.Hori. The first non-trivial case was
worked out by Seidel \cite{seidel-lefschetz} who constructed a Fukaya-type
$A_\infty$-category associated to a symplectic Lefschetz pencil. Since
this construction is a model for the general definition we briefly
review it next.

Let $(\underline{Y},\omega)$ be an open symplectic manifold, and let
$\bw : \underline{Y} \to \mathbb{C}$ be a symplectic Lefschetz
fibration, i.e.\ a $C^\infty$ complex-valued function with isolated
non-degenerate critical points $p_1,\dots,p_r$. This means that the
smooth parts of the fibers of $\bw$ are symplectic submanifolds of
$(\underline{Y},\omega)$, and that near each $p_{i}$ we can find a
$\omega$-adapted almost complex structure on $\underline{Y}$ so that
in almost-complex local coordinates $\bw$ is given by
$\bw(z_1,\dots,z_n)=\bw(p_i)+z_1^2+ \dots+z_n^2$.  Fix a regular value
$\lambda_0$ of $\bw$, and consider an arc $\gamma\subset \mathbb{C}$
joining $\lambda_0$ to a critical value $\lambda_i=\bw(p_i)$.  Using
the horizontal distribution which is symplectic orthogonal to the
fibers of $f$, we can transport a cycle vanishing at $p_i$ along the
arc $\gamma$ to obtain a Lagrangian disc $D_\gamma\subset
\underline{Y}$ fibered above $\gamma$, whose boundary is an embedded
Lagrangian sphere $L_\gamma$ in the fiber
$\underline{Y}_0=\bw^{-1}(\lambda_0)$. The Lagrangian disc $D_\gamma$
is called the {\it Lefschetz thimble} over $\gamma$, and its boundary
$L_\gamma$ is the Lagrangian vanishing cycle associated to the
critical point $p_i$ and to the arc $\gamma$.

Let $\gamma_1,\dots,\gamma_r$ be a collection of arcs in $\mathbb{C}$ joining
the reference point $\lambda_0$ to the various critical values of $\bw$,
intersecting each other only at $\lambda_0$, and ordered in the
clockwise direction around $p_0$. Each arc $\gamma_i$ gives rise to a
Lefschetz thimble $D_i\subset \underline{Y}$, whose boundary is a Lagrangian
sphere $L_i\subset \underline{Y}_0$. After a small perturbation we can always
assume that these spheres intersect each other transversely inside
$\underline{Y}_0$. Following Seidel \cite{seidel-lefschetz} we have

\begin{defi}\label{def:fs} Fix a commutative unital
  coefficient ring $R$.  The directed Fukaya-Seidel category $\op{{\sf
FS}}\left((\underline{Y},\omega,\bw);\{\gamma_i\}\right)$ over $R$ is
the following $A_\infty$-category:
\begin{itemize}
\item the objects of
  $\FS\left((\underline{Y},\omega,\bw);\{\gamma_i\}\right)$ 
are the Lagrangian vanishing cycles $L_1,\dots,L_r$; 
\item the morphisms
between the objects are given by
$$\mathrm{Hom}(L_i,L_j)=\begin{cases}
CF^*(L_i,L_j;R)=R^{|L_i\cap L_j|} & \mathrm{if}\ i<j\\
R\cdot id & \mathrm{if}\ i=j\\
0 & \mathrm{if}\ i>j;
\end{cases}$$
and the differential $m_1$, composition $m_2$ and higher order
products $m_k$ are defined in terms of Lagrangian-Floer homology inside
$\underline{Y}_0$. 
\end{itemize} 

More precisely,
$$m_k:\mathrm{Hom}(L_{i_0},L_{i_1})\otimes \dots\otimes
\mathrm{Hom}(L_{i_{k-1}},L_{i_k}) \to \mathrm{Hom}(L_{i_0},L_{i_k})[2-k]$$
is trivial when the inequality $i_0<i_1<\dots<i_k$ fails to hold (i.e.\ it
is always zero in this case, except for $m_2$ where composition with an
identity morphism is given by the obvious formula).

When $i_0<\dots<i_k$, $m_k$ is defined by fixing a generic
$\omega$-compatible almost-complex structure on $\underline{Y}_0$ and counting 
pseudo-holomorphic maps from a disc with $k+1$ cyclically ordered
marked points on its boundary to $\underline{Y}_0$, mapping the marked points
to the given intersection points between vanishing cycles, and the portions
of boundary between them to $L_{i_0},\ldots,L_{i_k}$ respectively.

The derived Fukaya-Seidel category
$\FS\left((\underline{Y},\omega,\bw);\{\gamma_i\}\right)$ is defined as
the homotopy category of the category of  twisted complexes (see
\cite{bondal-kapranov} for the definition) over $\op{{\sf
FS}}\left((\underline{Y},\omega,\bw);\{\gamma_i\}\right)$.
\end{defi}

\

\bigskip

Next we turn to the holomorphic ($B$-brane) side of the picture and
discuss the $B$-twisted Landau-Ginzburg model. Following \cite{O} we define

\begin{defi} Let $Y$ be a complex algebraic variety and let
  $\bw : Y \to \mathbb{C}$ be a regular function. The derived
category $D^b(Y,\bw)$ of the potential $\bw$ is defined as the
disjoint union of the quotient categories $Q(Y_{c})$ of all singular
fibers $Y_{c} := \bw^{-1}(c)$, $c \in \op{crit}(\bw)$. The category
$Q(Y_{c})$ is the quotient category of derived category of coherent
sheaves $D^b(Y_{c})$ modulo the full triangulated subcategory of
perfect complexes on $Y_{c}$.
\end{defi}

The categories $D^b(Y,\bw)$ and
$\FS\left(\underline{Y},\omega,\bw\right)$ together with the derived
category $D^b(X)$ of coherent sheaves of a complex variety $X$ and the
Fukaya category $\Fuk(\underline{X},\omega)$ play a central role in
the HMS conjecture for non-Calabi-Yau manifolds. In particular, when
the symplectic manifold $(\underline{X},\omega)$ underlies a
canonically polarized variety of general type and
$(\underline{X},\omega)\mir (Y,\bw)$ or $X\mir
(\underline{Y},\omega,\bw)$ are mirror pairs, the HMS conjecture
predicts a precise relationship between $\Fuk(\underline{X},\omega)$
and $D^b(Y,\bw)$, as well as between $D^{b}(X)$ and
$\FS\left(\underline{Y},\omega,\bw\right)$. Next we investigate this
relationship in more detail.

Suppose $(\underline{X},\omega)$ is a compact symplectic manifold
which underlies some complete intersection in a toric variety.  The
standard procedure for constructing the Landau-Ginzburg partner of
$(\underline{X},\omega)$ is based on the Hori-Vafa ansatz
\cite{HV}. For every such $(\underline{X},\omega)$ the Hori-Vafa
prescription produces an affine  complex Landau-Ginzburg mirror
$(Y^{\op{aff}},\bw^{\op{aff}})$ where $Y^{\op{aff}}$ is a closed
algebraic submanifold in an affine torus $(\mathbb{C}^{*})^{n}$, and
the superpotential $\bw^{\op{aff}}$ is the restriction of a regular
function on $(\mathbb{C}^{*})^{n}$. Physically the properties of
$(\underline{X},\omega)$ are encoded in the critical level sets of the
mirror superpotential $\bw^{\op{aff}}$. Often the affine mirror data
$(Y^{\op{aff}},\bw^{\op{aff}})$ is too crude as an invariant and does
not capture all the information in the $A$-model on
$(\underline{X},\omega)$. This inadequacy of
$(Y^{\op{aff}},\bw^{\op{aff}})$ usually
manifests itself geometrically in the fact that the critical loci of
$\bw^{\op{aff}}$ are not compact and to obtain a viable Landau-Ginzburg
theory one has to partially compactify $Y^{\op{aff}}$ so that $\bw^{\op{aff}}$
becomes a proper map \cite[Section~7.3]{HV}.  Suppose $\bw^{\op{prop}} :
Y^{\op{prop}} \to \mathbb{C}$ is such partial compactification.
Usually $Y^{\op{prop}}$ is a variety with complicated
singularities. The next step is to find a crepant resolution $Y$ of
$Y^{\op{prop}}$. After doing this and extending appropriately the
superpotential, one gets a projective morphism $\bw: Y \to \mathbb{C}$
with total space $Y$ which is a smooth non-compact Calabi-Yau.

The pair $(Y,\bw)$ is the Landau-Ginzburg partner of
$(\underline{X},\omega)$. By definition the category $D^b(Y,\bw)$ of
$B$-branes for $(Y,\bw)$ depends only on the structure of the singular
fibers of $\bw$. The Landau-Ginzburg partner we construct is not
uniquely defined by any means - but this was not to be expected
either. Even in the case of Calabi-Yau manifolds the mirror partners
are not uniquely defined. In our setup, there is also the additional
issue of the different ways of building $Y$. Our numerous computations
show that no matter how one (partially) compactifies and
desingularizes the affine Hori-Vafa mirror, at least combinatorially
the structure of the critical levels is the same. We expect that this
will be the case on the level of categories as well. That is, we
expect that $D^b(Y,\bw)$ will depend on $(\underline{X},\omega)$
only. When $X$ is 3-dimensional this is indeed true: the possible
desingularizations of $Y^{\op{prop}}$ differ by flops, and the
category we get is the same (see Section~\ref{sec:flops}).

In the following example we illustrate our approach for constructing a
Landau-Ginzburg partner, and show that the produced categories have
the expected by HMSC behaviour. The example deals with a del Pezzo
surface in which case HMSC is a theorem \cite{AKO1}. We have chosen
this example since the degree four del Pezzo is an intersection of two
quadrics in $\mathbb{P}^{4}$ and the analysis of its mirror
illustrates well how mirror constructions for intersection of quadrics
should be approached.

\begin{ex}
Following the procedure from \cite{HV} (see also the begininning of
our section \ref{sec-main}) we compute the mirror of the manifold
$(\underline{X},\omega)$ which underlies a degree four Del Pezzo
surface $X$ realized as the intersection of two quadrics in
$\CP^{4}$. We get the following affine Hori-Vafa miror:
$$
Y^{\op{aff}} =
\begin{cases}
x_{1}\cdot x_2\cdot x_3\cdot x_4\cdot x_5  = e^{-t} =: A \\
x_1 + x_2  = -1 \\
x_ 3 + x_{4}  = -1 
\end{cases}$$
and a potential 
\[
\bw^{\op{aff}} = x_{5}.
\]
If we solve for $x_{2}$ and $x_{4}$ from the second and third equation
and for $x_{5}$ from the first equation, we get 
\[
Y^{\op{aff}} = \mathbb{C}^{*}\times \mathbb{C}^{*}, \text{ with coordinates
  $x_{1}$ and $x_{3}$} 
\]
and 
$$
\bw^{\op{aff}} = \frac{A}{x_1(x_1+1)x_3(x_3+1)}.
$$ 
In order to partially compactify and resolve the total space of the
potential we do the following. Since the compatification is defined
modulo birational transformations, we can start with any
compactification of $\mathbb{C}^{*}\times \mathbb{C}^{*}$. We will use
the four-point blowup of $\CP^1 \times \CP^1$ where the points of the
blow ups are $(x_1,x_3)=(\infty,0), (\infty,-1), (0,\infty),
(-1,\infty)$.  On this blow-up the potential becomes a morphism to
$\CP^{1}$ and so we will get a fiberwise compactification $\bw : Y \to
\mathbb{C}$ if we delete the fiber at infinity, i.e. the proper
transforms of the four lines $x_{1} = 0$, $x_{3} = 0$, $x_{1} = -1$,
and $x_{3} = -1$ where the morphism $\CP^{1}\times \CP^{1} \to
\CP^{1}$ has poles.  This gives us a superpotential whose critical
locus consists of one isolated point $(x_1,x_3)=(-1/2,-1/2)$, and a
degenerate critical locus consisting of the proper transforms of the
lines $x_1=\infty$ and $x_3=\infty$.

We make a change of coordinates $t=1/x_1$ and $u=1/x_3$ , so that $\bw =
t^2 u^2 / (t+1)(u+1)$ and the critical points are $(-2,-2)$ and the
proper transforms of the two coordinate axes $t=0$ and $u=0$, while
the poles are at $t$ or $ u = -1$ or $\infty$.  We blow up $(0,-1)$, $
(0,\infty)$, $(-1,0)$, $(\infty,0)$.  We can deform $\bw$ to a new
superpotential $\bw' = t(t-e)u(u-e)/(t+1)(u+1)$, where $e$ is a small
constant. This deformation preserves the compactness property -
namely, $\bw'$ is well-defined over the blown-up space where we were
defining $\bw$, and its critical points lie inside a compact subset, and
do not escape to infinity when $e$ tends to $0$.  In fact, one can
also check that the topology of the generic fiber isn't changed (in
all cases it is a four-times punctured elliptic curve).  The potential
$\bw'$ has only isolated critical points -- 8 of them. Therefore we can
compute the vanishing cycles for $\bw'$. As usual we do
not count critical points over the fiber over $\infty$ in $ \CP^1 $.

The 8 critical points of $\bw'$  are:

\begin{enumerate}
\item[1)] in the 4-component fiber $\bw'$ =0, the four points $(0,0)$,
  $(0,e)$, $(e,0)$, $(e,e)$.
\item[2)] for non-zero values of $\bw'$, the critical points correspond
to solving:
$$
\begin{cases}d \log \bw'/dt = 1/t + 1/(t-e) - 1/(t+1) = 0,
\\
d \log \bw'/du = 1/u + 1/(u-e) - 1/(u+1) = 0,
\end{cases}\ \mathrm{i.e.}\ \begin{cases}t^2 + 2t - e = 0\\
u^2 + 2u - e = 0\end{cases}
$$ 
which gives four points lying close to $(e/2,e/2)$, $(e/2,-2)$,
$(-2,e/2)$, $(-2,-2)$ (if $e$ is small). Note that $\bw'(e/2,e/2)\simeq
e^4/16$, $\bw'(e/2,-2) = \bw'(-2,e/2)\simeq e^2$, $\bw'(-2,-2)\simeq 16$.
\end{enumerate}

\smallskip

\noindent
So, we have $\bw'$ with 4 critical values
($\lambda_1=\dots=\lambda_4=0$, $\lambda_5\simeq e^4/16$,
$\lambda_6=\lambda_7\simeq e^2$, $\lambda_8\simeq 16$).  We fix a
reference fiber $\bw' = \lambda_0 = e^4/32$ , and arcs $\gamma_i$
joining $\lambda_0$ to the critical values, to define vanishing
cycles:

\

\begin{enumerate}
\item[1.]  For $L_1,\ldots , L_4$  we choose a straight line from
$\lambda_0$  to $0$. 
\item[2.] For $ L_5$      we choose a straight line from $\lambda_0$
to $\lambda_5$.
\item[3.] For $L_6, L_7$ we choose an arc from $\lambda_0$ to
$\lambda_6=\lambda_7$ passing below the real axis.
\item[4.] For $L_8$ we choose an arc from $\lambda_0$ to $\lambda_8$
passing below the real axis.
\end{enumerate}

\

We can find the vanishing cycles by viewing each fiber (except for
$\bw'=0$ which is too singular) as a double cover of the $t$-plane
branched at 4 points.  The branch points of the fiber $\bw = \lambda$
are given by roots of: $\lambda^2 (t+1)^2 + e^2 t^2 (t-e)^2 + (4+2e)
\lambda t (t-e) (t+1) = 0$, which we can plot for various values of
$\lambda$ to find the vanishing cycles.  We will just summarize the
answer:

\begin{itemize}
\item $L_1,\ldots , L_4$  are  all disjoint (they occur in the same
  fiber $\bw'=0$). 
\item $L_1,\ldots,L_4$ all intersect each of $L_5,\ldots , L_8$ in
  exactly one point. 
\item $L_5$  intersects $L_6$  and $L_7$  in 2 points, $ L_8$  in 4
  points. 
\item $L_6$  and $ L_7$  don't intersect; $ L_8$  intersects $ L_6$
  and $ L_7$  in 2 points. 
\end{itemize}

\

The  matrix of the corresponding quiver is:
\[
\begin{pmatrix}
1 & 0 &  0 &  0 &  1 &  1 & 1 & 1 \\
  & 1 &  0 &  0 &  1 &  1 & 1 & 1 \\
  &   &  1 &  0 &  1 &  1 & 1 & 1 \\
  &   &    &  1 &  1 &  1 & 1 & 1 \\
  &   &    &    &  1 &  2 & 2 & 4 \\
  &   &    &    &    &  1 & 0 & 2 \\
  &   &    &    &    &    & 1 & 2 \\
  &   &    &    &    &    &   & 1
\end{pmatrix}
\]
All entries correspond to actual numbers of morphisms, and the answer
perfectly matches the derived category of a 5-point blowup of $\CP^2$
(or equivalently a 4-point blowup of $\CP^1\times\CP^1$), which is
indeed a Del Pezzo surface of degree $4$.

According to definition \ref{def:fs} the Fukaya-Seidel category is the
$A_{\infty}$-category of modules over the $A_{\infty}$-version of the
path algebra of this quiver where the higher products are given by
disk instantons.  With additional work, one can check that $m_1=0$ and
$m_k=0$ $\forall k\ge 3$, and compute the composition $m_2$. So we
have the following \cite{AKO1}:

\begin{theo} The categories $D^b(X)$ and
  $\FS(\underline{Y},\bw, \eta)$ are equivalent. 
The same statement holds true for the derived categories of
noncommutative deformations of $X$ and the categories associated with
deformed symplectic structures on $\bw : \underline{Y} \to \mathbb{C}$.
\end{theo}
\end{ex}

\section{Main conjectures} \label{sec-main}

Let $X$ be a manifold of
general type (i.e. $X$ has a positive-definite canonical class).
In this section we will 
formulate and motivate an analogue of the HMS conjecture for such manifolds.  

A quantum sigma-model with target $X$ is free in the infrared limit,
while in the ultraviolet limit it is strongly coupled.  In order to
make sense of this theory at arbitrarily high energy scales, one has
to embed it into some asymptotically free $N=2$ field theory, for
example into a gauged linear sigma-model (GLSM).  Here ``embedding''
means finding a GLSM such that the low-energy physics of one of its
vacua is described by the sigma-model with target $X$. In mathematical
terms, this means that $X$ has to be realized as a complete
intersection in a toric variety. 

The GLSM usually has other vacua as well, whose physics is not related
to $X$. Typically, these extra vacua have a mass gap. To learn about
$X$ by studying the GLSM, it is important to be able to recognize the
extra vacua. Let $Z$ be a toric variety defined as a symplectic
quotient of $\mathbb{C}^N$ by a linear action of the gauge group $G\simeq
U(1)^k$. The weights of this action will be denoted $Q_{ia}$, where
$i=1,\ldots,N$ and $a=1,\ldots,k$. Let $X$ be a complete intersection
in $X$ given by homogeneous equations $G_\alpha(X)=0$,
$\alpha=1,\ldots, m$. The weights of $G_\alpha$ under the G-action
will be denoted $d_{\alpha a}$. The GLSM corresponding to $X$ involves
chiral fields $\Phi_i,$ $i=1,\ldots,N$ and $\Psi_\alpha$,
$\alpha=1,\ldots,m$.  Their charges under the gauge group $G$ are
given by matrices $Q_{ia}$ and $d_{\alpha a}$, respectively.  The
Lagrangian of the GLSM depends also on complex parameters $t_a$,
$a=1,\ldots,k$. On the classical level, the vector $t_a$ is the level
of the symplectic quotient, and thus parametrizes the complexified
K\"ahler form on $Z$.  The K\"ahler form on $X$ is the induced one. On
the quantum level the parameters $t_a$ are renormalized and satisfy
linear RG equations:
$$
\mu\frac{\partial t_a}{\partial\mu}=\beta_a=\sum_i Q_{ia}.
$$
In the Calabi-Yau case all $\beta_a$ vanish, and the parameters $t_a$ are not renormalized. 

The mirror Landau-Ginzburg model has (twisted) chiral fields $\Lambda_a$,
$a=1,\ldots,k,$ $Y_i$, $i=1,\ldots,N$ and $\Upsilon_\alpha,$
$\alpha=1,\ldots,m$. The superpotential is given by
$$ 
\bw=\sum_a \Lambda_a\left(\sum_i Q_{ia}Y_i-\sum_\alpha d_{\alpha a}
\Upsilon_\alpha - t_a\right) + \sum_i e^{-Y_i} +\sum_\alpha
e^{-\Upsilon_\alpha}.
$$

The vacua are in one-to-one correspondence with the critical points of
$\bw$. By definition, massive vacua are those corresponding to
non-degenerate critical points. An additional complication is that
before computing the critical points one has to partially compactify
the target space of the Landau-Ginzburg model.

One can determine which vacua are ``extra'' (i.e. unrelated to $X$) as
follows. The infrared limit is the limit $\mu\to 0$. Since $t_a$
depend on $\mu$, so do the critical points of $\bw$. A critical point
is relevant for $X$ (i.e. is not an extra vacuum) if and only if the
critical values of $e^{-Y_i}$ all go to zero as $\mu$ goes to zero.
In terms of the original variables $\Phi_i$, this means that vacuum
expectation values of $|\Phi_i|^2$ go to $+\infty$ in the infrared
limit. This is precisely the condition which justifies the classical
treatment of vacua in the GLSM. It is instructive here to recall that
that the classical space of vacua in the GLSM is precisely $X$.

Now let us state the analogue of the HMS for complete intersections
$X$ which are of general type.

\begin{conn} \label{con-AtoB} 
If $(\underline{X},\omega)|(Y,\bw)$ is a Hori-Vafa mirror
  pair, then:
\begin{itemize} 
\item The Fukaya category
$\Fuk(\underline{X},\omega)$ embeds as a direct summand into the
category $D^b(Y,\bw)$ of $B$-branes of $\bw : Y \to
\mathbb{C}$.
\item The orthogonal complement of
$\Fuk(\underline{X},\omega)$ in $D^b(Y,\bw)$ is very simple: it is a
direct sum of several copies of the category of graded modules over a
Clifford algebra of a symmetric bilinear form on a complex vector
space of dimension $n=\dim_{\mathbb{C}} Y$.
\end{itemize}
\end{conn}

\

\noindent
There is also an analogous mirror conjecture in which the $A$ and $B$
sides of the theory are switched. More precisely if denote the
$C^{\infty}$-manifold underlying $Y$ by $\underline{Y}$, then the
complex manifold $X$ of general type should determine a symplectic
structure $\eta$ on $\underline{Y}$ with respect to which $\bw :
\underline{Y} \to \mathbb{C}$ becomes a symplectic Lefschetz fibration.

\

\noindent
With this notation we have the following:

\begin{conn} \label{con-BtoA} 
If $X|(\underline{Y},\bw,\eta)$ is such a mirror pair, then the
category $D^{b}(X)$ of $B$-branes on $X$ is equivalent to the
$A$-brane category $\FS(\underline{Y}_{D},\bw,\eta)$ of a potential
$\bw : \underline{Y}_{D} \to D$, where $0 \in D \subset \mathbb{C}$ is
a suitably chosen disk, and $\underline{Y}_{D} = \bw^{-1}(D)$.
\end{conn}                           

\

Let us explain the physical reasoning behind these conjectures. To any
$N=(2,2)$ quantum field theory in two dimensions one can assign the
categories of $A$-branes and $B$-branes.  These categories are
``topological'', in the sense that they do not depend on the flat
two-dimensional metric necessary to defined the quantum field
theory. Since rescaling the metric is equivalent to rescaling the
renormalization scale $\mu$, this means that the categories of
topological D-branes are $\mu$-independent. By construction, the GLSM
and its mirror Landau-Ginzburg model have identical physics in the
infrared limit $\mu\to 0$, up to a mirror involution, so their
categories of topological D-branes must be equivalent, up to exchange
of $A$ and $B$-branes.

It is important to note at this point that categories of topological
D-branes are well-defined even when the topological twisting is
not. For example, if $c_1(X)\neq 0$, then the $B$-type twist for the
sigma-model with target $X$ does not exist, but the category of
$B$-branes on $X$ is well-defined.  The reason for this is that the
structure of the category is specified by the boundary chiral ring,
whose definition uses only the properties of the field theory on a
flat world-sheet (with boundaries).  If the topological twist exists,
it provides additional structures on the category of branes (such as a
Serre functor).

The neighborhood of each vacuum can be regarded as a separate physical
theory, with its own categories of topological D-branes, which are
full sub-categories in the respective categories of D-branes in the
``parent'' theory.  This implies that the category of $A$-branes on
$X$ is a full sub-category of the category of $B$-branes for the
mirror Landau-Ginzburg model, and vice versa.

The category of $B$-branes for a Landau-Ginzburg model has very strong
localization properties~\cite{orlov-completion}: it is a direct sum of
categories determined by ``critical-level schemes'' of $\bw$. In a
generic situation the values of $\bw$ at all extra vacua are different
and also different from the values of $\bw$ at relevant vacua. This
implies that each extra vacuum contributes independently, and that the
total category of $B$-branes for the Landau-Ginzburg model is a sum of
a category equivalent to the category of $A$-branes on $X$ and the
contributions of extra vacua.

Finally, if a critical point of $\bw$ is non-degenerate, then the
corresponding category of $B$-branes is equivalent to the category of
graded modules over the Clifford algebra with $n$ generators, where
$n$ is the number of variables on which $\bw$ depends. We summarize the
discussion above in the Table~\ref{table:HMSC} below.

\begin{table}[ht!]
\begin{center}
\begin{tabular}[t]{|c|c|}
\hline A side & B side \\ \hline\hline  
\begin{minipage}[c]{\miniwidth}
\

\medskip

$\underline{X}$ --- compact manifold, \\
$\omega$ --- symplectic form on $\underline{X}$.

\

\medskip

\end{minipage}
&
\begin{minipage}[c]{\miniwidth}
$X$ --- smooth projective variety over $\mathbb{C}$
\end{minipage} 
\\ 
\begin{minipage}[c]{\miniwidth} 
\[\Fuk(X,\omega)=
\left\{
  \begin{array}{l}
\tn{Obj: } (L_i,\shE)\\ \tn{Mor: } HF(L_i,L_j)
  \end{array}
\right.
\]
$L_i$ --- Lagrangian submanifold of $X$,\\ $\shE$ --- flat $U(1)$-bundle on
$L_i$ 

\

\medskip

\end{minipage}
&
\begin{minipage}[c]{\miniwidth}
\[D^b(X)=\left\{
  \begin{array}{l}
\tn{Obj: } C_i^\bullet\\ \tn{Mor: } Ext(C_j^\bullet,C_j^\bullet)
  \end{array}
\right.
\]
$C_i^\bullet$ --- complex of coherent sheaves on $X$
\end{minipage}
\\ \hline \multicolumn{2}{c}{ \xymatrix{
\ar@{<->}[drr]&&\ar@{<->}[dll]\\ && } } \\ \hline 
\begin{minipage}[t]{\miniwidth}

\

\medskip

$(\underline{Y},\eta)$ --- open symplectic manifold, \\
${\bw :\underline{Y}\to\mathbb{C}}$ ---
a  proper $C^{\infty}$ map with symplectic fibers. 
\end{minipage}
&
\begin{minipage}[t]{\miniwidth}

\

\medskip

$Y$ --- smooth quasi-projective variety over $\mathbb{C}$, \\
$\bw :Y\to \mathbb{C}$ ---
  proper algebraic map. 
\end{minipage}
\\
\begin{minipage}[t]{\miniwidth}
\[
\FS\left(\underline{Y},\bw,\eta\right)=\left\{
  \begin{array}{l}
\tn{Obj: } (L_i,\shE)\\ \tn{Mor: } HF(L_i,L_j)
  \end{array}
\right.
\]
$L_i$ --- Lagrangian submanifold of $\underline{Y}_{\lambda_{0}}$,\\ 
$\shE$ --- flat $U(1)$-bundle on $L_i$. 
\end{minipage}
&
\begin{minipage}[t]{\miniwidth}
\[D^b(Y,\bw)=\bigsqcup_t D^b(Y_t)\Big/\op{Perf}(Y_t)
\]
\end{minipage}
\\
\begin{minipage}[c]{\miniwidth}
$\FS_{\lambda_{i},r_{i}}(\underline{Y},\bw,\eta)$ -- the
  Fukaya-Seidel category of $(\underline{Y}_{|t - \lambda_{i}| <
  r_{i}},\bw,\eta)$.
\end{minipage}
&
\begin{minipage}[c]{\miniwidth}
\[
D^b_{\lambda_{i},r_{i}}(Y,\bw)= \bigsqcup_{|t -
  \lambda_{i}|\le r_{i}}
D^b(Y_t) \Big/ \op{Perf}(Y_t)
\]

\

\medskip
\end{minipage}
\\ \hline 
\end{tabular}
\caption{Kontsevich's HMS conjecture} \label{table:HMSC}
\end{center}
\end{table}

In this formalism we need to take a Karoubi closure on both sides of
HMSC.

\

\bigskip

\noindent
The above statement of HMS interacts in a subtle way with various
other physical and geometric operations, such as orientifold
projections, phase flows in gauged linear sigma models, or large
$N$-duality. In this direction
we focus on a conjecture relating the complex or symplectic geometry of
a pair of manifolds, one of which is a Fano, and the other is of
general type:

\begin{conn} \label{con:side.A.equiv} \ Let $F$ and $G$ be
  projective manifolds. Suppose $F$ is a Fano manifold, and $G$ is a
  manifold of general type or a Calabi-Yau manifold. Suppose further
  that there exists a fully faithful functor 
\[
\Phi : D^{b}(G) \hookrightarrow
  D^{b}(F)
\] 
between the derived categories of $F$ and $G$. Write
  $\underline{F}$ and $\underline{G}$ for the $C^{\infty}$ manifolds
  underlyng $F$ and $G$. Then:
\begin{itemize}
\item[(i)] For every complexified K\"{a}hler class $\alpha_{F}$ on
  $F$ there exists a complexified  K\"{a}hler class $\alpha_{G}$ on
  $G$ and a fully faithful functor 
\[
\Psi : \overline{\Fuk}(\underline{G},\alpha_{G}) \hookrightarrow
\overline{\Fuk}(\underline{F},\alpha_{F})
\]
between the corresponding Karoubi completed Fukaya categories.
\item[(ii)] If $K_{\Phi} \in D^{b}(G\times F)$ is a kernel object for $\Phi$
  (such an object exists by \cite{orlov-K3}), then $\Psi$ is given by a kernel
  object $K_{\Psi}$in $\overline{\Fuk}(\underline{G}\times
  \underline{F},-p_{G}^{*}\alpha_{G} + p_{F}^{*}\alpha_{F})$ which is
  uniquely determined by $K$.
\end{itemize}
\end{conn}

\

\noindent
It is hard to make this conjecture more precise mainly because at the
moment we do not have a good enough grasp on all objects in the Fukaya
category of a symplectic manifold. Heuristically, the kernel object
$K_{\Psi}$ should be a coisotropic $A$-brane in the sense of
\cite{ko-coisotropic}, and one should be able to write the connection
data for this coisotropic brane in terms of a (super) connection data
on $K_{\Phi}$ which satisfies an appropriate form of the
Hermite-Yang-Mills equations and is uniquely determined by the
holomorphic structure on $K_{\Phi}$ and the K\"{a}hler structure
$-p_{G}^{*}\alpha_{G} + p_{F}^{*}\alpha_{F})$ on $\underline{G}\times
\underline{F}$.


\section{Hyperelliptic curves} \label{sec-hyperelliptic}

In this section we study various aspects of Conjectures~\ref{con-AtoB}
and \ref{con-BtoA} in the case of hyperelliptic curves. In particular
we give two explicit models for the Landau-Ginzburg mirror of a genus
two curve.  

\subsection{Curves of genus two} \label{ssec-genus2}

Let $\underline{C}$ be a differentiable compact Riemann surface of
genus two. One presentation of such a surface is as the $C^{\infty}$
manifold underlying a divisor of bi-degree $(2, 3)$ in the Hirzebruch
surface ${\mathbb F}_0 = \mathbb{P}^1\times \mathbb{P}^1$. In other words
$\underline{C}$ is defined as the zeros of a generic form
$p(z_1:z_2;w_1:w_2)\,\, \in \,\,\Gamma ({\mathbb F}_0, \mathcal{O}(2,
3))$. Here $z_1, z_2, w_1$, and $w_2$ are the homogeneous coordinates
of the rulings of ${\mathbb F}_0$. If we realize the latter as a GIT 
quotient ${\mathbb C}^4 /\!/ {\mathbb C}^{\times 2}$, then the weights of
the torus action are given in the table:

\begin{center}
\begin{tabular}{|c|c|c|c|c|    } \hline
$z_{1}$ & $z_{2}$ & $w_{1}$ & $w_{2}$ & $p$ \\ \hline \hline
$1$ & $1$ &  $0$ & $0$ & $-3$ \\ \hline
$0$ & $0$ &  $1$ & $1$ & $-2$ \\ \hline

\end{tabular}
\end{center} 

\

\noindent
Let $B + i \omega$ be a complexified K\"{a}hler form on the
toric surface $\mathbb{F}_{0}$, which is homogeneous under the natural
torus action. In other words $B$ is a torus invariant closed real two form on
$\mathbb{F}_{0}$, and $\omega$ is a torus invariant 
K\"{a}hler form.  The Hori-Vafa
mirror of $\left(\underline{C},B + i \omega\right)$ is
constructed in stages:
\begin{itemize}
\item[(i)] Build the affine Hori-Vafa mirror $\bw^{\op{aff}} :
  Y^{\op{aff}} \to \mathbb{C}$ of $\left(\underline{C},B
  + i \omega\right)$ corresponding to the embedding in
  $\mathbb{F}_{0}$; 
\item[(ii)] Compactify partially
  $\left(Y^{\op{aff}},\bw^{\op{aff}}\right)$ to a family
  $\bw^{\op{prop}} : Y^{\op{prop}} \to
  \mathbb{C}$, so that $Y^{\op{prop}}$ is a quasi-projective
  Gorenstein variety with trivial canonical class, and
  $\bw^{\op{prop}}$ is proper;
\item[(iii)] Construct a crepant (possibly stacky) resolution $\bw : Y
  \to \mathbb{C}$ of the compactified family
  $\bw^{\op{prop}} : Y^{\op{prop}} \to \mathbb{C}$;
\end{itemize} 
In fact from the point of view of Conjectures~\ref{con-AtoB} and
\ref{con-BtoA}  these construction steps should be supplemented by 
\begin{itemize}
\item[(iv)] Select a disk $0 \in D \subset \mathbb{C}$ so that the restricted
  family $(Y_{D} := \bw^{-1}(D),\bw)$ is the mirror of the genus two curve.
\end{itemize}
Explicitly this means that:
\begin{description}
\item[($A$ to $B$)] There is an equivalence \label{page:AtoB} of
  Karoubi closures
\[
\overline{\Fuk}(\underline{C},B+i\omega) \cong
  \overline{D^{b}}(Y_{D},\bw) = \coprod_{t \in D}
  \overline{D^{b}(Y_{t})/\op{Perf}(Y_{t})}; 
\]
\item[($B$ to $A$)] The differentiable manifolds and functions
  underlying the holomorphic Landau-Ginzburg mirrors for the
  various $B + i\omega$ are all naturally diffeomorphic to a fixed
  $C^{\infty}$ fibration $\bw : \underline{Y} \to
  \mathbb{C}$. Moreover fpr any choice of a complex structure on
  $\underline{C}$, i.e. for any choice of a section
  $p(z_1:z_2;w_1:w_2)\,\, \in \,\,\Gamma ({\mathbb F}_0,
  \mathcal{O}(2,3))$, the complex curve $C_{p} : \{ p = 0\}$ has an
  associated complexified symplectic form $\eta_{p}$ on
  $\underline{Y}_{D} := \bw^{-1}(D)$ so that $\op{Im} \eta_{p}$
  restricts to a symplectic form on the fibers of $\bw$ and we
  have an equivalence
\[
D^{b}(C_{p}) \cong \overline{\FS}(\underline{Y}_{D},\bw,\eta_{p}).
\]
where the bar on the right hand side denotes the Karoubi closure.
\end{description}

\

\noindent
The first step in the constuction is just an implementation of the
Hori-Vafa recipe described in the previous section.  Following
\cite{HV} we define the affine mirror of $\left(\underline{C},B + i
\omega\right)$ to be the pencil $\bw^{\op{aff}} : Y^{\op{aff}} \to
\mathbb{C}$ given by the function $\bw= x_1 + \cdots + x_5$ over the
threefold $Y^{\op{aff}} \subset ({\mathbb C}^{*})^{5}$ defined by the
equations: 
\[
Y^{\op{aff}} \ : \
\begin{cases}
x_1\cdot x_2  = a_1\cdot x_3^3 \\ 
x_4\cdot x_5  = a_2 \cdot x_3^ 2.
\end{cases}
\]
Here $x_{1}$, $x_{2}$, $x_{3}$, $x_{4}$, $x_{5}$ are the natural
coordinates on $({\mathbb C}^{*})^{5}$, and $a_{1}$ and $a_{2}$ are
non-zero complex constants corresponding to the $(B+i\omega)$-volume
of the two rulings of $\mathbb{F}_{0}$. Concretely we have
\[
a_{1} := \exp\left(-\int_{\{\op{pt}\}\times \mathbb{P}^{1}}
(B+i\omega)\right), \text{ and } a_{2} :=
\exp\left(-\int_{\mathbb{P}^{1}\times \{\op{pt}\}} (B+i\omega)\right).
\]
The most delicate and involved part of the above prescription are the
steps (ii) and (iii), and carrying those out will occupy the rest of
the section. 

For step (ii) we note that equivalently we can think of $Y^{\op{aff}}$
as the subvariety in $(\mathbb{C}^{*})^{4}\times \mathbb{C}$ with
coordinates $(x_{1},x_{2},x_{3},x_{4};w)$, defined by the equations
\[
Y^{\op{aff}} \ : \ \left\{ \begin{split}
x_4^2 - x_4(w - x_1 - x_2 - x_3) + a_{2} x_3^2 & = 0, \\
- x_1x_2 +  a_{1}x_3^3 & = 0.
\end{split}\right.
\]
In these terms the superpotential $\bw$ becomes simply the projection
on the factor $\mathbb{C}$ in the product $(\mathbb{C}^{*})^{4}\times
\mathbb{C}$. To obtain the partial compactification
$(Y^{\op{prop}},\bw^{\op{prop}})$ we embed $(\mathbb{C}^{*})^{4}$ as
the standard torus inside $\mathbb{P}^{4}$, and take $Y^{\op{prop}}$
to be the closure of $Y^{\op{aff}}$ in $\mathbb{P}^{4}\times
\mathbb{C}$, while $\bw^{\op{prop}}$ is again the projection on the
factor $\mathbb{C}$.

Proceeding with step (iii) we have to resolve the singularities of
$Y^{\op{prop}}$. The structure of the resolution is summarized in the
following:

\begin{theo} \label{theo-genus2}
There exists a crepant resolution $Y$ of $Y^{\op{prop}}$ so that the
the fibers of the function  $\bw : Y \to \mathbb{C}$
induced from $\bw^{\op{prop}}$  can be described as follows.
\begin{itemize}
\item For $w \neq 0, \frac{1}{a_{2}}\left(\frac{1\pm
  2\sqrt{a_1}}{3}\right)^3$,  the
fiber $Y_w := \bw^{-1}(w)$ is an elliptic K3 surface with two singular
fibers of type $I_9$ and six singular fibers of type $I_1$. So, $Y_w$
is a K3 surface of Picard rank 18.
\item For $w = \frac{1}{a_{2}}\left(\frac{1\pm
  2\sqrt{a_1}}{3}\right)^3$, the fiber $Y_w$ is singular and has
  a single isolated $A_{1}$-singularity. Its minimal resolution is an
  elliptic K3 with a $2\cdot I_{9} + 4\cdot I_{1} + I_{2}$ 
configuration of singular fibers.
\item The central fiber $Y_0$ of the fibration $\bw :Y \to \mathbb{C}$
is a union of three rational surfaces: $Y_0 = \widetilde{S}_1 \cup
\widetilde{S}_2 \cup \widetilde{S}_{3}$. These surfaces intersect in
three copies of $\mathbb{P}^{1}$:
\[
q_1 = \widetilde{S}_1 \cap \widetilde{S}_{3}, \qquad
q_2 = \widetilde{S}_2 \cap \widetilde{S}_{3}, \quad \text{and}
\quad q_3 = \widetilde{S}_1 \cap \widetilde{S}_2.
\] 
The three curves $q_{i}$ all meet in two points: $q_1
\cap q_2 \cap q_3 = \{M, N\}$. In particular $Y_0$ is
a type III degeneration of a K3 surface. The
geometry of the components of $Y_0$ and the position of the curves
$q_i$ on them is described as follows.
\begin{enumerate}
\item[(i)] The component $\widetilde{S}_{3}$ is a blow up of a
Hirzebruch surface $\mathbb{F}_2$ at 4 points of depth 2 on two
generic sections of $\mathbb{F}_2$ and such that two pairs of these 4
points belong to two (generic) fibers of $\mathbb{F}_2$. The proper
transforms of the two sections are the curves $q_1$ and
$q_2$. Their intersection points are $M$ and $N$.
\item[(ii)] To describe the union $\widetilde{S}_{1}\cup
\widetilde{S}_{2}$ we start as in (i) by blowing up $\mathbb{F}_{2}$
at the same points but only once this time. Consider the proper
transforms of the two sections as in (i) but call them here $q'$ and
$q''$. Take two copies of thus blown up $\mathbb{F}_2$ and identify
them along $q'$. The curves $q''$ will be identified with
$q_1$ and $q_2$ on $F$. To get the union
$\widetilde{S}_{1} \cup \widetilde{S}_{2}$ we have to blow up two more
times either the first or the second $\mathbb{F}_2$ at points
infinitesimally near the already blown up ones on $q'$.\footnote{The
fact that this construction of $\widetilde{S}_{1}\cup
\widetilde{S}_{2}$ is legitimate follows from the fact that (before
blowing the two last points up) there is an obvious $\ZZ_2$ symmetry
in the preimage of $Y_0^{\op{prop}}$ which exchanges
$\widetilde{S}_{1}$ and $\widetilde{S}_{2}$ and leaves the component
$\widetilde{S}_{3}$ invariant. This in turn shows that the choices of
the sections on the $\mathbb{F}_2$ we are making have to be in certain
agreement with each other.}
\end{enumerate}
 The rank of the Picard group of $Y_0$ is 21.
\end{itemize}
\end{theo}
{\bfseries Proof:} As explained above we can obtain the partial
compactification $\bw^{\op{prop}} : Y^{\op{prop}} \to \mathbb{C}$ by
closing $Y^{\op{aff}}$ in $\mathbb{P}^{4}\times
\mathbb{C}$. By slightly abusing notation we will write
$(x_{0}:x_{1}:x_{2}:x_{3}:x_{4})$ for the the homogeneous coordinates
on $\mathbb{P}^{4}$ and $\bw$  for the coordinate on $\mathbb{C}$. In
these terms we have
\[
Y^{\op{prop}} \ : \ \left\{
\begin{split}
x_4^2 - x_4(\bw x_{0} - x_{1} - x_{2} - x_{3}) + a_{1}x_3^2 & = 0 \\
- x_{0}x_{1}x_{2} +  a_{2}x_3^3 & = 0 
\end{split}
\right.
\]
and $\bw^{\op{prop}}$ is induced by the projection to the second factor of
$(\C^*)^4\times \C$. 

First we need to understand the general fiber of the desingularization
of $Y^{\op{prop}}$:

\begin{claim}
There is a crepant resolution 
\[
\xymatrix@C-1pc{
Y \ar[rr] \ar[rd]_-{\bw} & & Y^{\op{prop}} \ar[ld]^-{\bw^{\op{prop}}}
\\
& \mathbb{C} &
}
\]
of the singularities of $Y^{\op{prop}}$ so that:
\begin{itemize}
\item 
For $w \neq 0,\; \frac{1}{a_{2}}\left(\frac{1\pm
  2\sqrt{a_1}}{3}\right)^3$, the fiber $Y_w = \bw^{-1}(w)$ is an elliptic
K3 surface with two singular fibers of type $I_9$ and 6 singular fibers
of type $I_1$. 
\item For $w = \frac{1}{a_{2}}\left(\frac{1\pm
  2\sqrt{a_1}}{3}\right)^3$, the fiber $Y_w$ has an isolated
  singularity o type $A_{1}$ and its minimal resolution is an elliptic
  K3 with a $2\cdot I_{9} + 4\cdot I_{1} + I_{2}$ configuration
  ofsingular fibers.
\end{itemize}
\end{claim}
{\bfseries Proof.}
The singular locus of $Y^{\op{prop}}$ consists of nine
sections of $\bw^{\op{prop}}$: three of 
transversal singularity type $A_1$
\[
l_1: (1: w: 0: 0: 0; w) \qquad l_2: (1: 0: w: 0: 0; w)
\qquad  l_3: (0: 1: -1: 0: 0; w),
\]
and six of transversal type $A_2$
\[
\begin{array}{lll}
m_1:  (1: 0: 0: 0: w; w) & m_2: (1: 0: 0: 0: 0; w)
&  m_3: (0: 0: 1: 0: 0; w) \\
m_4:  (0: 1: 0: 0: 0; w) & m_5: (0: 0: 1: 0: -1; w)
&  m_6: (0: 1: 0: 0: -1; w).
\end{array}
\]
The singularities of $Y^{\op{prop}}$
can be resolved by blowing up these lines in the ambient space 
$\mathbb{P}^{4}\times \mathbb{C}$.

The fiber $Y^{\op{prop}}_w$ is transversal to the components of the
singular locus of $Y^{\op{prop}}$ for $w \neq 0$, and hense inherits
$3A_1$ and $6A_2$ singular points which form the singular locus of the
generic fiber of $\varphi$. These points get desingularized when we
blow up the lines. From this it folows that the critical points of $\bw :
(Y - {Y_0}) \to (\C - {0})$ come from the ones of $\bw^{\op{prop}} :
Y^{\op{prop}} \to \C$. A direct computation shows that the latter map
has only two critical points, at $w = \frac{1}{a_{2}}\left(\frac{1\pm
  2\sqrt{a_1}}{3}\right)^3$.

Before turning to the structure of the central fiber $Y_0$, let us note
one more feature of the manifold $Y$. If we denote by $S$ the surface
$\{ x_0x_1x_2 = a_{2}x_3^3 \} \subset \CP^3$, then $Y$ is
a double cover of 
$S\times \C$ via the map which forgets the variable $x_4$ (that is,
via projection with center $(0: 0: 0: 1) \in \CP^4$
fiberwise). Accordingly, every fiber $Y_w$ is a double cover of
$S$. This additional structure of $Y$ will help us understand the
geometry of the fibers $Y_w$.

Consider first $w \neq 0$. The cubic surface $S$ has three $A_2$ 
singularities which are the pairwise intersection points of the three
lines, 
\[
t_0: x_0 = x_3 = 0,  \qquad t_1: x_1 = x_3 = 0, \qquad t_2: x_2 = x_3.
= 0.
\]
Note that these lines are all contained $S$ and that the triangle
formed by this lines is the intersection of $S$ with the plane $x_{3}
= 0$. The fiber $Y^{\op{prop}}_{w}$ is a double cover of $S$ branched
at the union of two cubic curves on $S$:
\[
D_1 = \{w\cdot x_0 - x_1 - x_2 - x_3 - 2\sqrt{a_{1}}x_{3} = 0\} \cap
S, \qquad 
D_2 = \{w\cdot x_0 - x_1 - x_2 - x_3 + 2\sqrt{a_{1}}x_{3} = 0\} \cap S
\]  
which intersect transversally at the three points 
\[
P_0 = (0: 1: -1: 0), \qquad  P_1(1: 0: w: 0), \qquad  P_2(1: w: 0: 0)
\]
belonging to $t_0, t_1$, and $t_2$ respectively. In particular, we see
that there are two sources of singularities for $Y^{\op{prop}}_{w}$:
the singular locus of $S$, and the points of $D_1 \cup D_2$. Hence
$Y_w$ is a double cover of the surface $S'$ obtained by first
de-singularizing $S$ and then blowing up additionally the points $P_0,
P_1, P_2$. The surface $S'$ is a rational elliptic surface. It has one
fiber of type $I_9$ and three fibers of type $I_1$. For $w =
\frac{1}{a_{2}}\left(\frac{1 + 2\sqrt{a_1}}{3}\right)^3$, one of these
$I_1$ fibers is $D_1$, while for $w = \frac{1}{a_{2}}\left(\frac{1-
2\sqrt{a_1}}{3}\right)^3$ one of them is $D_2$. The surface $Y_w$ is a
double cover of $S'$ branched at the fibers of $S'$ corresponding to
$D_1$ and $D_2$.  As a result we get that $Y_w$ is an elliptic
fibration with two $I_9$ and six $I_1$ fibers for $w \neq
\frac{1}{a_{2}}\left(\frac{1\pm 2\sqrt{a_1}}{3}\right)^3$, while for
$w = \frac{1}{a_{2}}\left(\frac{1\pm 2\sqrt{a_1}}{3}\right)^3$ the
double cover has a singular point of type $A_1$ at the preimage of the
node of $D_{1}$ and $D_{2}$ respectively. The double cover
interpretation also shows that the generic $Y_w$ is a smooth K3
surface of Picard number 18. This proves the claim. \ \hfill $\Box$

\

\bigskip

We turn now to the case $w = 0$ and the description of the central
fiber. The fiber $Y^{\op{prop}}_{0}$ has six singular points: $4A_2$
and $1A_1$ singular points at infinity ($x_0 = 0$), and a deeper
singular point at $(1: 0: 0: 0: 0; 0)$ which corresponds to the
intersection point of the components $l_1, l_2, m_1,$ and $m_2$ of the
singular locus of $Y^{\op{prop}}$. In terms of the double cover
$Y^{\op{prop}}_{0} \to S$ we see that the branch divisor has two
singular points: $P_0$ (at infinity) and the intersection $t_1 \cap
t_2$. The latter point corresponds to a collision of $P_1$ and $P_2$
as the generic branch divisor specializes to the the branch divisor of
$Y^{\op{prop}}_{0} \to S$. As before we see that the points at
infinity get resolved with the blowing up the lines $l_3, m_3, m_4,
m_5,$ and $m_6$. The new feature here is the resolution at the point
$(1: 0: 0: 0: 0; 0)$.

Consider $Y^{\op{prop}}$ away from $x_0 = 0$. We will denote this
threefold by $V^a := Y^{\op{prop}} - \{ x_{0} = 0\} $. It 
is a toric variety which in the $\C^5$ with coordinates $x_{1}$,
$x_{2}$, $x_{3}$, $x_{4}$, $x_{5}$ is given by the equations
\[
V^{a} \ : \  \left\{ \begin{split} 
x_1x_2 & = a_{1}x_3^3\\ 
x_4x_5 & = a_{2}x_3^2.
\end{split} \right.
\] 
Note also that the lines $l_1, l_2, m_1$, and $m_2$ (which all meet at
the origin) are toric subvarieties of $V^a$. The desingularization of
$Y^{\op{prop}}$ at the origin can be easily understood in terms of the toric
geometry of $V^a$. Indeed, the Newton polytope of $V^a$ is shown in
the picture below.

\

\begin{figure}[!ht]
\begin{center}
\epsfig{file=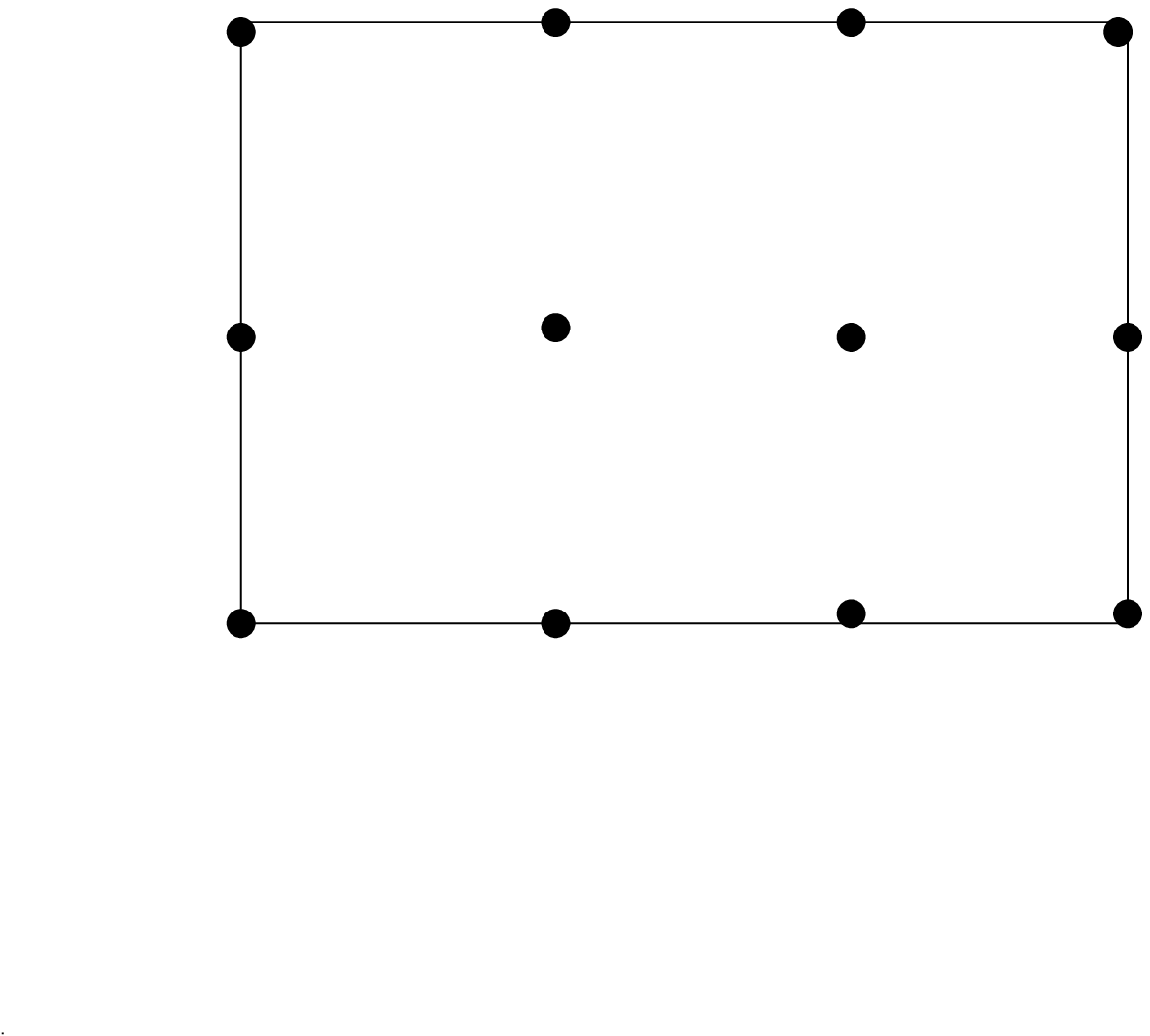,width=2in} 
\end{center}
\vspace{-4pc}
\caption{The Newton polytope of $V^a$.}
\label{fig-M1} 
\end{figure}

\

\noindent
Desingularization of $V^a$ amounts to subdividing this polytope. There
are different ways to subdivide the rectangle. One can get from one
subdivision to another geometrically by perfoming flops on the
desingularized varieties. We choose to subdivide the rectangle by
first blowing up $m_2$. Then we blow up (in any order) the proper
preimages of the lines $m_1, l_1$, and $l_2$. The result after these
blow ups can be explained as follows.

\bigskip

\

\begin{figure}[!ht]
\begin{center}
\epsfig{file=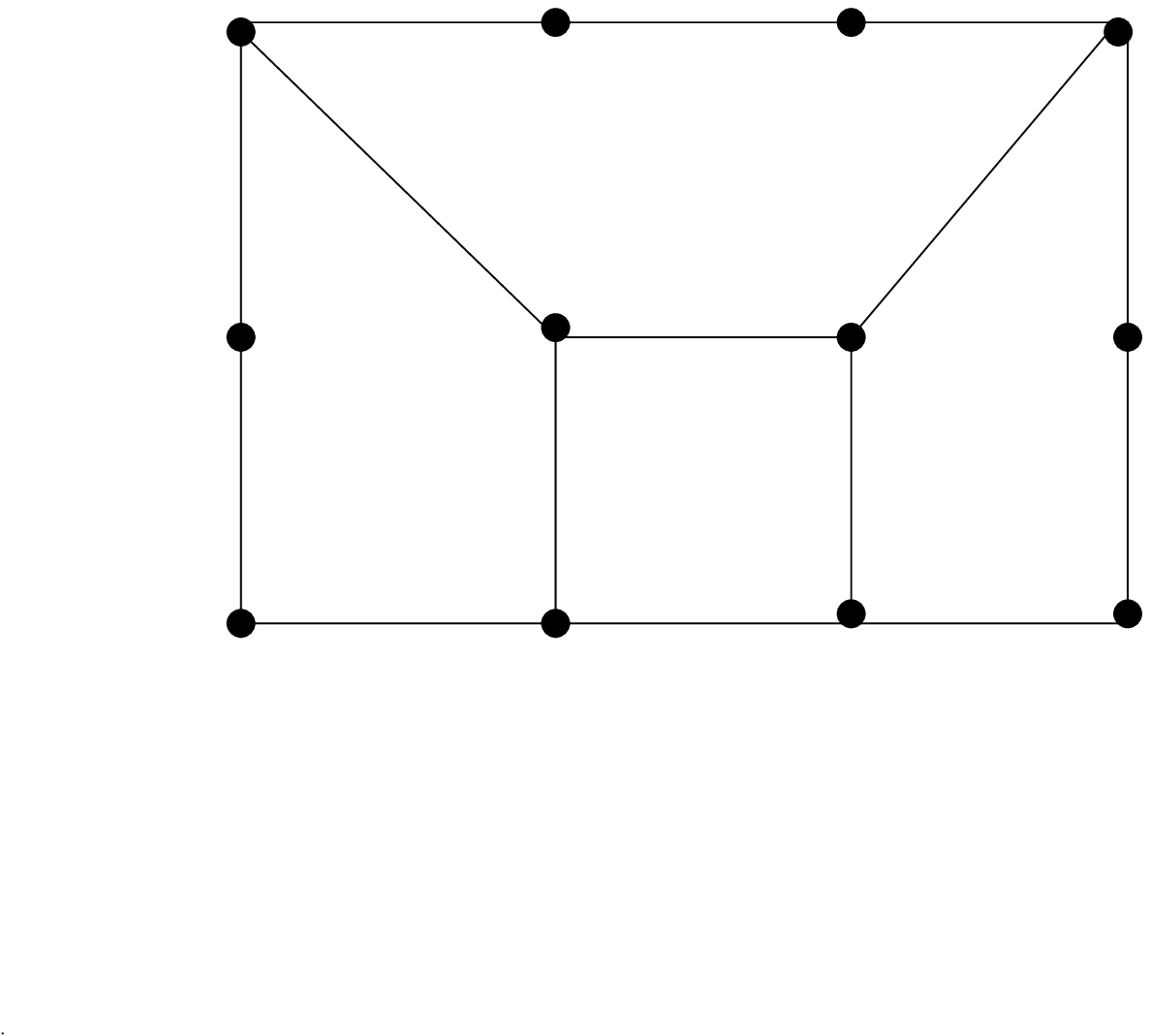,width=2in} 
\end{center}
\vspace{-4pc}
\caption{The Newton polytope of $V^b$.}
\label{fig-M2} 
\end{figure}

\noindent
(i) Denote by $V^b$ the proper transform of $V^a$ after we blow up
$m_2$ in $\mathbb{C}^{5}$. The Newton polytope of $V^{b}$ is depicted in
Figure~\ref{fig-M2}.

\

The central fiber $V_0^b$ is a union of three components: two
exceptional $\CP^2$'s, call them $H_1$ and $H_2$, and the proper
transform $F^a$ of $V^a_0$. Going back to the double cover structure
on $Y^{\op{prop}}_0$, we see that $V^a_0$ is a double cover of $S^a =
S - \{x_0 = x_3 = 0\} \subset \C^3$ branched at two curves which
intersect at the origin of $\C^3$. These curves intersect in such a
way that the blow up of $m_2$ takes them appart, and $F^a$ is actually
a double cover of a smooth surface $\widetilde{S}^a$ (the only
singularity of $S^a$ gets resolved) branched in a smooth divisor. We
see then that the components of $V^b_0$ are all smooth surfaces.

We are now ready to describe the proper transform $F$ of
$Y^{\op{prop}}_{0}$ that we get after we resolve the whole
$Y^{\op{prop}}$. As mentioned above, $Y^{\op{prop}}_0$ is a double
cover of $S$ branched in two curves, and the singularities of
$Y^{\op{prop}}_0$ arise from the singularities of $S$ and the
singularities of the branch curve. On the surface $S$ the efect of
blowing up the lines $m_1, \, m_3, \, m_4,\, m_5,\,$ and $\,m_6$ in
the ambient space, can be seen as a resolution of the singularities of
the branch curve.  Call the resulting surface $\widetilde{S}'$. The
proper transform of $Y^{\op{prop}}_0$ under these modifications is a
double cover of $\widetilde{S}'$ branched at two curves which
intersect at a point which belongs to the line $l_3$. After resolving
along $l_3$, we obtain the ultimate proper transform $F$ as a double
cover of the blow up of $\widetilde{S}'$ at the point where $l_3$
intersects it, branched in a smooth divisor. So $F$ is itself
smooth. The surface $\widetilde{S}'$ can also be obtained from $\CP^2$
by blowing up two of the vertices of a triangle three times each so
that the side these vertices define gets blown up only twice, and the
other two get blown up three times each. The branch curves are then
the preimages of two lines, in general position w.r.t. the triangle,
and with a commom point on the side of the triangle blown up only
twice. Accordingly, after we blow up along $l_3$, we can identify
$\widetilde{S}'$ with the Hirzebruch surface $\mathbb{F}_1$ blown up
six times at points belonging to a fiber and two sections, in general
position, so that the sections get blown up three times each. From
this point of view the branch divisors will consist of two general
fibers of the blown up $\mathbb{F}_{1}$. As a result we get a
description of $F$as the blowing up of two general sections on a
Hirzebruch surface $\mathbb{F}_2$ at points belonging to two general
fibers (four points in all) and their infinitesimally near ones so
that the sections get blown up three times in each of the ponts. We
still need to resolve along $m_1,\, l_1,\,$ and $l_2$. As before, the
toric picture will help us understand how this resolution will affect
the picture we already have. We will complete this in (ii) and (iii)
below.

\bigskip

\noindent
(ii) Next we study  how $F$, $H_1$ and $H_2$ intersect each other.
Consider the blow up of $\C^{5}$ at $m_2$ in the affine chart where 
\[
x_1 = u_1/u_4, \qquad x_2 = u_2/u_4, \qquad x_3 = u_3/u_4.
\]
The equations of the proper transform of $V^a$ are
\[
V^{b} \ : \ \left\{ 
\begin{split}
 (u_1/u_4)\cdot(u_2/u_4) = a_{1}\cdot x_4\cdot(u_3/u_4)^3\\ w -
x_4\cdot(1 + u_1/u_4 + u_2/u_4 + u_3/u_4) = a_{2}\cdot
x_4\cdot(u_3/u_4)^2
  \end{split} \right. 
\]
In $Y$ the exceptional divisor of this blow up is given by $x_4
= w = (u_1/u_4)\cdot(u_2/u_4) = 0$, so it is a union of two
planes which we called $H_1$ and $H_2$ respectively.  In the $\CP^3$
with homogeneous coordinates $(u_1: \dots: u_4)$. 
The proper transform $V^b$ has a
fiber above $w=0$ consisting of three components: the exceptional
planes $H_1,\,H_2\,$, and the proper transform of $V^a_0$.

The planes $H_1$ and $H_2$ intersect in the line $q_3 : \; u_1 = u_2 =
0 \subset \CP^3$. The intersection of
$H_1$ (respectively  $H_2$) with $F$ is  a curve $q_1$ (respectively $q_2$)
which is the preimage of one of the exceptional curves of the
desingularization $\widetilde{S}$ of $S$. These $(-2)$ exceptional curves,
as obvious from the explanations above, meet the branch divisor of the
double cover $F \to \widetilde{S}$ at two points each. This is why, on
$F$, they will lift to $\CP^1$'s of self-intersection $-4$. The curve
$q_i$ is a conic on $H_i$, $\,\,i = 1, 2$. All three curves intersect
in two triple points.

For future references let us point out that $H_i$ considered as toric
subvarieties of $V^b$ have three fixed points and three fixed lines
each under the torus action. Two of the points are common (belong to
$q_3$, which is one of the fixed lines for each of the planes), and
the respective third points belong to the individul planes. As it
follows by direct inspection, the conics $q_1$ and $q_2$ intersect
$q_3$ away from the fixed points on the latter, and the former lines
pass through the third fixed point on the corresponding
plane. Finally, the lines $u_4 = 0$ on each of $H_1$ and $H_2$ are
tangent at the torus-fixed points to $q_1$ and $q_2$ while the
analogous lines with $u_3 = 0$ intersect the conics in two points (one
of them is a torus-fixed one). Also, as it can be computed in the
coordinates above, the $(-1)$-curves on $\widetilde{S}$ that result from
the blow-up of $\mathbb{F}_1$ lift to two pairs of $(-1)$-curves on $F$
in such a way that they intersect also the lines with $u_3 = 0$ in
each  $H_i$.

\bigskip

\noindent
(iii) Finally we explain the desingularization of $V^b$. The proper
transforms of $m_1, l_1$ and $l_2$ intersect $V^b_0 = H_1 \cup H_2 \cup
F$ in singular points of $V^b_0$. Hence, they intersect it in points
of the curves $q_i, \,\, i = 1, 2, 3$. These proper transforms are
toric subvarieties of $V^b$, so they intersect $H_i$ in torus-fixed
points. By the toric picture of $V^b$ we see that $H_1$ and $H_2$
intersect only one of the lines $l_1$ and $l_2$, and both intersect
the line $m_1$. Hence, $l_1$ and $l_2$ intersect the exceptional
planes at the torus-fixed points on the conics $q_i$, while $m_1$
intesects $q_3$ at the fixed point with $u_4 \neq 0$.

\begin{figure}[!ht]
\begin{center}
\epsfig{file=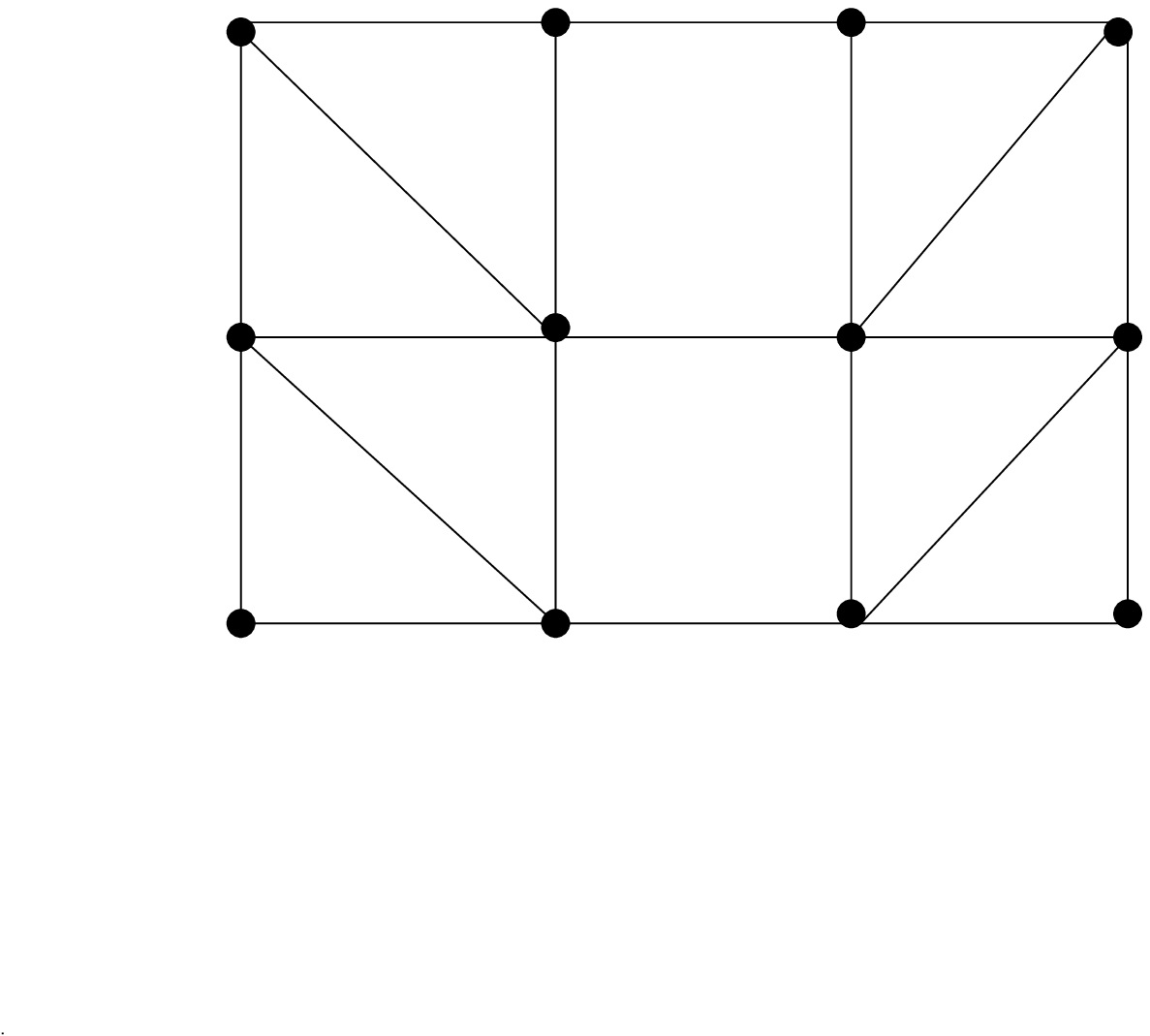,width=2in} 
\end{center}
\vspace{-4pc}
\caption{The Newton polytope of the resolution of $V^b$.}
\label{fig-M3} 
\end{figure}

Blowing up the three lines (and thus resolving the generic fiber of
$Y^{\op{prop}}$) blows up $H_1$ and $H_2$ in two points each (they
become del Pezzo surfaces, $\widetilde{H}_{1}$ and
$\widetilde{H}_{2}$, of degree 7 each). The common line of these del
Pezzos, the proper transform of $q_3$, has self-intersection $0$ on
each of them. In this process $F$ gets blown up in two points (to the
surface $F^b$), and the curves of intersection of this surface with
$\widetilde {H}_{1}$ and $\widetilde{H}_{2}$ have self-intersection
$3$ on the del Pezzos and $-5$ on $F^b$. The subdivision of the Newton
polytope we get at this point are shown on Figure~\ref{fig-M3}.

\begin{figure}[!ht]
\begin{center}
\epsfig{file=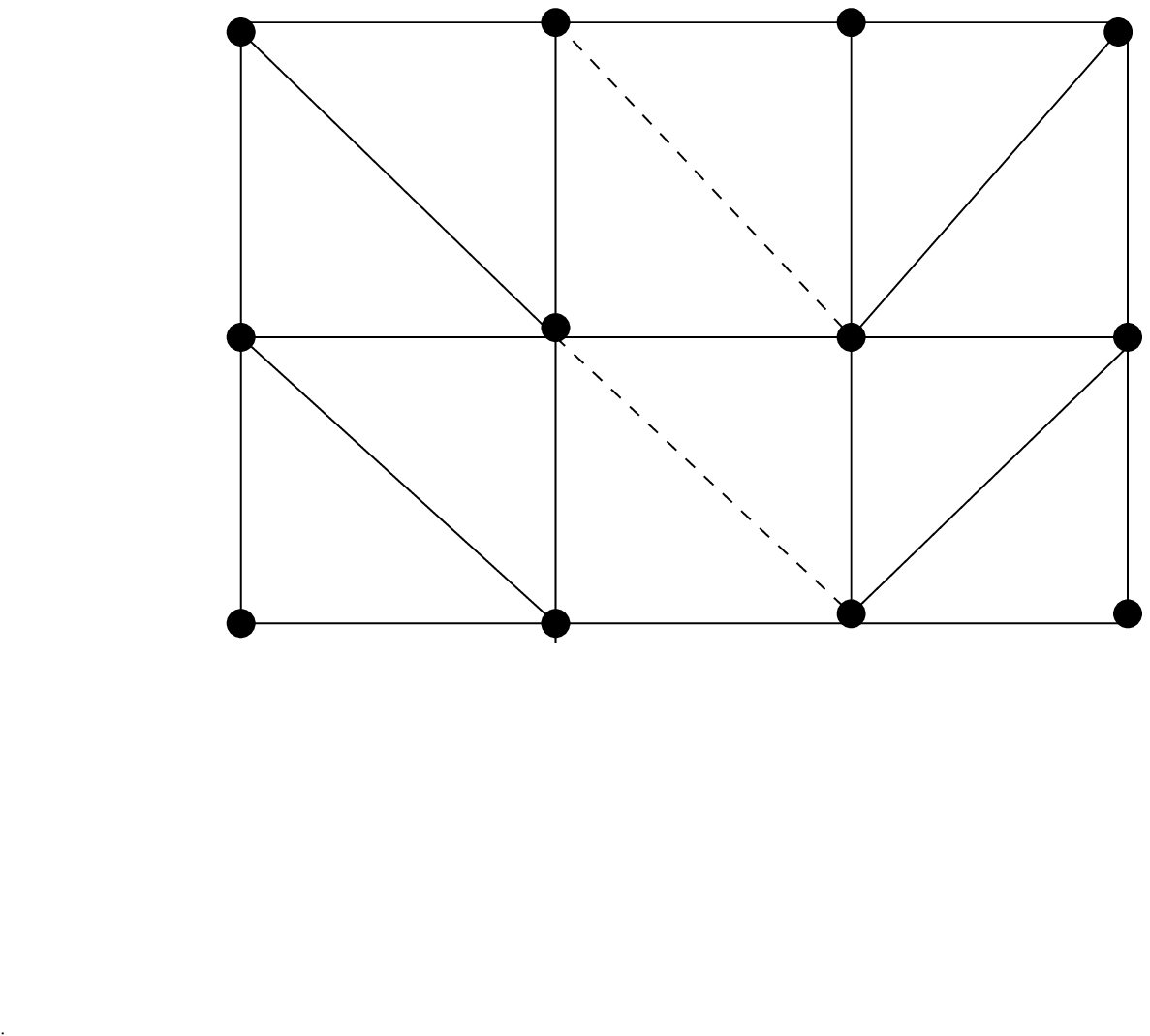,width=2in} 
\end{center}
\vspace{-4pc}
\caption{The Newton polytope of the small resolution.}
\label{fig-M4} 
\end{figure}

There are still two singular points on $V^b$ to be resolved. According
to the toric picture these are the common points of
$\widetilde{H}_{1}$ and $\widetilde{H}_{2}$ only. These actually are
the fixed points under the torus action on the intersection of the two
degree 7 del Pezzos (and are far from $F^b$). We see also that the two
points can be resolved by using small resolutions as
Figure~\ref{fig-M4} illustrates.

The resulting intersection line of the two exceptional components,
call the components $S_1$ and $S_2$, will be either $0, -2$ or $ -1,
-1$ or $-2, 0$ depending on the choice of the small resolution. Do now
two flops \label{page:flops} in the exceptional curves on $S_3$ from
the last blow up. The proper preimage of $Y^{\op{prop}}_0$ will then
be the surface $F$, and the exceptional surfaces $S_1$ and $S_2$ get
blown up one more time. The central fiber $Y_{0}$ of the
desingularization $Y$ will consist of three components:
$\widetilde{S}_{3}$, $\widetilde{S}_{1}$ and
$\widetilde{S}_{2}$. These intersect in three curves
($\mathbb{P}^{1}$s) with self-intersections $(-4, 2)$ for $(F,
\widetilde{S}_{i})$, and $(-2, 0)$, $(-1, -1)$, or $(0, -2)$ on
$(\widetilde{S}_{1}, \widetilde{S}_{2})$.

The component $\widetilde{S}_{3}$ corresponding $F$, is then a blow up
of the Hirzebruch surface $\mathbb{F}_2$ in the following sequence of
14 points. Choose two generic sections of $\mathbb{F}_2$ and mark on
them four points defined by the intersection with two generic
fibers. The proper transform $\widetilde{S}_{3}$ of
$Y^{\op{prop}}_{0}$ is the blow up of $\mathbb{F}_2$ at each of the
marked points, then at infinitesimally near points of order one and
two to the marked ones where, at each step, the infinitesimally near
points belong to the proper transforms of the chosen sections. The
proper transforms of the chosen sections after these 12 blow ups are
the curves of intersection $q_1$ and $q_2$ of $F$ with the exceptional
components $H_1$ and $H_2$. The component $\widetilde{S}_3$ is the
blow up of $F$ at two more points (the effect on $F$ of blowing up the
lines $l_1$ and $l_2$ on $V$). These are points on $q_1$ and $q_2$
which belong to exceptional curves of the blown up $\mathbb{F}_2$ as
explained above.

For the sake of later considerations, we prefer to work with a
different desingularization of $Y^{\op{prop}}$. As noticed above,
different desingularizations differ by flops from the one we just
constructed. Flops can be performed along $(-1, -1)$-lines, and these
are exactly the exceptional curves in $Y_0$ which meet only two
components of $Y_0$. In particular, we can do flops 4 times on $F$ in
such a way that this component gets obtained from $\mathbb{F}_2$ by
blowing up the two sections in points of depth two instead of
three. This will result in blowing up $\widetilde{S}_i$ two more times
(at points on $q_i$). For the sake of symmetry we also prefer to do
the small resolutions discussed above in both components $S_1$ and
$S_2$. We will be considering this desingularization and be using the
same notation for the central fiber and its components.  

\ \hfill $\Box$

\

\bigskip

\noindent As a corollary, we get the combinatorial structure of the
singularities of $\bw : Y \to \mathbb{C}$.

\begin{figure}[!ht]
\begin{center}
\psfrag{P}[c][c][1][0]{{$M$}}
\psfrag{Q}[c][c][1][0]{{$N$}}
\psfrag{Q1}[c][c][1][0]{{$q_{1}$}}
\psfrag{Q2}[c][c][1][0]{{$q_{2}$}}
\psfrag{Q3}[c][c][1][0]{{$q_{3}$}}
\psfrag{t1}[c][c][1][0]{{$y_{+}$}}
\psfrag{t2}[c][c][1][0]{{$y_{-}$}}
\epsfig{file=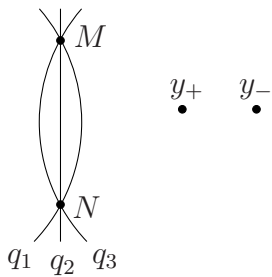,width=1.3in} 
\end{center}
\caption{The critical points of $\bw$.}
\label{fig-wW1} 
\end{figure}

\begin{cor}  \label{cor:genus.two}
The set of critical points of $\bw$ consists of:
\begin{itemize}
\item two isolated points
$y_{\pm}$ with $\bw(y_{\pm}) = \frac{1}{a_{2}}\left(\frac{1 \pm
2\sqrt{a_1}}{3}\right)^3$, and 
\item three $\mathbb{P}^1$'s: $q_1$, $q_2$,
$q_3$, passing through two other points $M$ and $N$. (See
figure~\ref{fig-wW1})). 
\end{itemize}
Note that the positive dimensional part of the critical locus,
i.e. the union $q_{1}\cup q_{2}\cup q_{3}$, is
a singular curve of 
genus 2.
\end{cor}

\

\begin{rem} \label{rem:karoubi.genus.two}
$\bullet$ \ From the proof of Theorem~\ref{theo-genus2} we see that
  $y_{\pm}$ are ordinary double points in their respective fibers and
  so the categories $D^{b}_{\op{sing}}(Y_{\bw(y_{\pm})}$ are just
  categories of modules over a Clifford algebra. Therefore for step
  (iv) above it is natural to take the disk $D \subset \mathbb{C}$ to
  be any disk which contains $0 \in \mathbb{C}$ bud does not contain
  the critical values $\frac{1}{a_{2}}\left(\frac{1 \pm
  2\sqrt{a_1}}{3}\right)^3$

\

\noindent
$\bullet$ \ In view of the completion theorem
\cite[Theorem~2.10]{orlov-completion} the Karoubi closure
$\overline{D^b}(Y,\bw)$ of the category $D^b(Y,\bw)$ is determined by
the points $y_{\pm}$ and a formal neighborhood of the singular curve
$q_1 \cup q_2 \cup q_3$ inside the fiber $Y_{0}$. 
\end{rem}

\

\noindent
Our next task is to analyze the the structure of the category
$D^{b}(Y_{D},\bw) = D^{b}_{\op{sing}}(Y_{0})$ in more detail and in
particular to understand how this category depends on the moduli of
the Landau-Ginzburg model $(Y,\bw)$.  To do that we will need a better
model for $(Y,\bw)$ in which the components of the singular fiber
$Y_{0}$ appear in a more symmetric way. 

Since $Y$ is three dimensional, the category of a superpotential $\bw
: Y \to \mathbb{C}$ does not change if we modify $Y$ by a simple flop
(see Section~\ref{sec:flops}) 
we can try to simplify the zero fiber of our superpotential by
performing flops on the compactification constructed in
Theorem~\ref{theo-genus2}. Before we can explain the outcome of these
modifications we will need the following

\

\medskip

\noindent {\bfseries Construction.} \ Fix a complex number
$\lambda$. Consider the ruled surface $\mathbb{F}_{2}$ with its
natural projection $\pi : \mathbb{F}_{2} :=
\mathbb{P}(\mathcal{O}(2)\oplus \mathcal{O}) \to \mathbb{P}^{1}$.  Let
$\delta \subset \mathbb{F}_{2}$ be the unique $(-2)$-section of
$\pi$. Fix four distinct points $p_{0}$, $p_{\infty}$, $p_{1}$,
$p_{2}$ of cross ration $(p_{0},p_{\infty},p_{1},p_{2}) = \lambda$,
and let $P_{1}$, $P_{2}$ be points in $\mathbb{F}_{2} - \delta$ with
$\pi(P_{i}) = p_{i}$, $i = 1,2$. Consider the pencil of curves in the
linear system $|\mathcal{O}(\delta)\otimes \pi^{*}\mathcal{O}(2)|$
passing through $P_{1}$ and $P_{2}$. Let $Q'$ and $Q''$ be two
irreducible curves that generate this pencil. Note that $Q'$ and $Q''$
are sections of $\pi$ that do not intersect $\delta$ and so they each
have self intersection $2$. Note that if we fix the points
$p_{0},p_{\infty},p_{1},p_{2}$ the automorphisms of $\mathbb{F}_{2}$
that act along the fibers of $\pi$ and preserving $\delta$ will act
transitively on the data $(P_{1},P_{2},Q',Q'')$. 

Indeed these
automorphisms are just affine automorphims of the total space of
$\mathcal{O}(2)$ and so are generated by scaling along the fibers or
translations by sections of $\mathcal{O}(2)$. Since we can always find
a quadratic polynomial to take two prescribed values at $p_{1}$ and
$p_{2}$, we can use a translation to bring $P_{1}$ and $P_{2}$ to
points on the zero section of $\mathcal{O}(2)$. After this translation
$Q'$ and $Q''$ will be just two distinct sections of $\mathcal{O}(2)$ which
vanishes at $p_{1}$, $p_{2}$. Since up to scale there is a unique
non-zero section $\xi$ of $\mathcal{O}(2)$ which
vanishes at $p_{1}$, $p_{2}$ we can find  constants $c' \neq c''$ so that
$Q' = c'\cdot \xi$  and $Q'' = c''\cdot \xi$. Translating by
$(-c'\cdot \xi)$ and scaling by $(c''-c')^{-1}$ will keep the points
$P_{1}$, $P_{2}$ on the zero section but will move the sections $Q'$
and $Q''$ to the zero section and the section $\xi$
respectively. 

In particular, up to isomorphism the data $(P_{1},P_{2},Q',Q'')$ is
determined by the four points $p_{0},p_{\infty},p_{1},p_{2}$. Consider
the rational surface obtained by blowing up $\mathbb{F}_{2}$ at the
four points on $Q'$ and $Q''$ sitting over the points $p_{0}$ and
$p_{\infty}$. By what we just explained up to isomorphism this surface
is uniquely determined by the ordered qudruple of points
$p_{0},p_{\infty},p_{1},p_{2}$ in $\mathbb{P}^{1}$ Since the quadruple
of points has only one modulus - the cross ratio $\lambda$ - it
follows that up to isomorphism the rational surface is determined by
$\lambda$. For future reference we will denote this surface by
$\mathfrak{S}^{\lambda}$. The proper transforms of $Q'$ and
$Q''$ are two rational curves on $\mathfrak{S}^{\lambda}$ intersecting
at two points. We will denote these curves 
by $\mathfrak{q}'$ and $\mathfrak{q}''$. Notice that the surfaces 
$\mathfrak{S}^{\lambda}$ and $\mathfrak{S}^{1/\lambda}$ are naturally
isomorphic since we can pass from one to the other by switching the
order of the points $p_{1}$ and $p_{2}$.

\

\bigskip

\noindent
We are now ready to describe the better model for the Landau-Ginzburg
mirror of $(\underline{C},B + i\omega)$:

\begin{theo} \label{theo-genus2-symmetric}
There exists a crepant resolution $\bfY$ of $Y^{\op{prop}}$ so that the
the fibers of the function  $\bfw : \bfY \to \mathbb{C}$
induced from $\bw^{\op{prop}}$  can be described as follows.
\begin{itemize}
\item For $\bfw \neq 0, \frac{1}{a_{2}}\left(\frac{1\pm
  2\sqrt{a_1}}{3}\right)^3$, the fiber $\bfY_w :=
  \bfw^{-1}(w)$ is an elliptic K3 surface with two singular fibers of
  type $I_9$ and six singular fibers of type $I_1$,
  i.e. $\bfY_w$ is a K3 
  surface of Picard rank 18.
\item For $w = \frac{1}{a_{2}}\left(\frac{1\pm
  2\sqrt{a_1}}{3}\right)^3$, the fiber $\bfY_w$ is singular and has
  a single isolated $A_{1}$-singularity. Its minimal resolution is an
  elliptic K3 with a $2\cdot I_{9} + 4\cdot I_{1} + I_{2}$ 
configuration of singular fibers.
\item The central fiber $\bfY_0$ of the fibration $\bfw :
\bfY \to \mathbb{C}$ is a union of three rational surfaces:
$\bfY_0 = \mathfrak{S}_1 \cup \mathfrak{S}_2 \cup
\mathfrak{S}_{3}$. These surfaces intersect in three copies of
$\mathbb{P}^{1}$:
\[
\bfq_1 = \bfS_1 \cap \bfS_{3}, \qquad
\bfq_2 = \bfS_2 \cap \bfS_{3}, \quad \text{and}
\quad \bfq_3 = \bfS_1 \cup \bfS_2.
\] 
The three curves $\bfq_{i}$ all meet in two points: $\bfq_1
\cap \bfq_2 \cap \bfq_3 = \{M, N\}$. Explicitly 
 $\bfY_0$ and the position of the curves
$\bfq_i$ on it can be described as follows. Fix once and for
all an ordering 
$\lambda_{1}$, $\lambda_{2}$ of the roots of the polynomial $x^{2} + x
+ a_{2} = 0$ and let $\lambda = \lambda_{1}/\lambda_{2}$ be their
ratio. Consider three copies $\bfS_{i}$, $i = 1,2,3$ of the
surface $\mathfrak{S}^{\lambda}$, each with two marked rational curves
$\mathfrak{q}'_{i}$ and $\mathfrak{q}''_{i}$. The fiber
$\bfY_{0}$ is the normal crossings surface obtained by gluing
the  components 
$\{ \bfS_{i}\}_{i = 1,2,3}$ along the curves
$\mathfrak{q}'_{i}$ and $\mathfrak{q}''_{i}$, where we identify 
$\mathfrak{q}''_{i}$ with $\mathfrak{q}'_{i+1}$  (for all $i \mod 3$)
and the identification are such that the ordered quadruples of marked
points on each curve are matched.

 The rank of the Picard group of $\bfY_0$ is 21.
\end{itemize}
\end{theo}
{\bf Proof.} The proof of the theorem follows immediately from the
proof of Theorem~\ref{theo-genus2}. We will not repeat any of the
arguments and will just indicate the flops that one needs in
order to construct $\bfY$ from $Y$. 

\begin{figure}[!ht]
\begin{center}
\psfrag{H1}[c][c][0.8][0]{{$H_{1}$}}
\psfrag{H2}[c][c][0.8][0]{{$H_{2}$}}
\psfrag{(1)}[c][c][0.8][0]{{$(1)$}}
\psfrag{(4)}[c][c][0.8][0]{{$(4)$}}
\psfrag{q1}[c][c][0.8][0]{{$q_{1}$}}
\psfrag{q2}[c][c][0.8][0]{{$q_{2}$}}
\psfrag{q3}[c][c][0.8][0]{{$q_{3}$}}
\psfrag{A}[c][c][0.8][0]{{$A$}}
\psfrag{B}[c][c][0.8][0]{{$B$}}
\psfrag{C}[c][c][0.8][0]{{$C$}}
\psfrag{D}[c][c][0.8][0]{{$D$}}
\psfrag{E}[c][c][0.8][0]{{$E$}}
\psfrag{F}[c][c][0.8][0]{{$F$}}
\psfrag{G}[c][c][0.8][0]{{$G$}}
\psfrag{H}[c][c][0.8][0]{{$H$}}
\epsfig{file=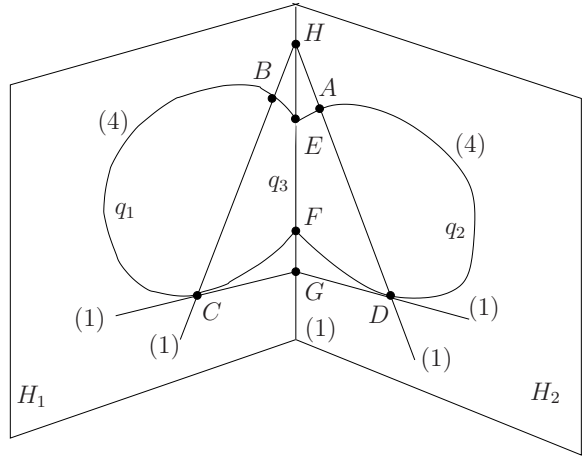,width=3in} 
\end{center}
\caption{The planes $H_{1}$ and $H_{2}$.}
\label{fig:H.with.conics} 
\end{figure}

After blowing up the line $m_{2}$ in $Y^{\op{prop}}$ two of the
components we had in the central fiber were the planes
$H_{1}$ and $H_{2}$ intersecting along the line $q_{3}$. Furthermore
each $H_{i}$ contained a conic $q_{i}$ as depicted on
Figure~\ref{fig:H.with.conics}.

\

\

\noindent
The self-intersection numbers of the various curves in
Figure~\ref{fig:H.with.conics} 
are indicated in parentheses. 
The points $H$, $G$, $C$, and $D$ appearing in
Figure~\ref{fig:H.with.conics} are fixed by the torus action. The
lines $GC$ and $GD$ are
tangent lines to the conics $q_{1}$ and $q_{2}$ at $C$ and $D$
respectively. The cross ratios $(D,F,E,A)$,
$(G,F,E,H)$, and $(C,F,E,B)$ are all equal to $\lambda$. This
observation will be essential in construction of $\bfY$. 

\begin{figure}[!ht]
\begin{center}
\psfrag{H1}[c][c][0.8][0]{{$H_{1}$}}
\psfrag{H2}[c][c][0.8][0]{{$H_{2}$}}
\psfrag{(-1)}[c][c][0.8][0]{{$(-1)$}}
\psfrag{(0)}[c][c][0.8][0]{{$(0)$}}
\psfrag{(2)}[c][c][0.8][0]{{$(2)$}}
\psfrag{(-2)}[c][c][0.8][0]{{$(-2)$}}
\psfrag{A}[c][c][0.8][0]{{$A$}}
\psfrag{B}[c][c][0.8][0]{{$B$}}
\psfrag{C}[c][c][0.8][0]{{$C$}}
\psfrag{D}[c][c][0.8][0]{{$D$}}
\psfrag{E}[c][c][0.8][0]{{$E$}}
\psfrag{F}[c][c][0.8][0]{{$F$}}
\psfrag{G}[c][c][0.8][0]{{$G$}}
\psfrag{H}[c][c][0.8][0]{{$H$}}
\epsfig{file=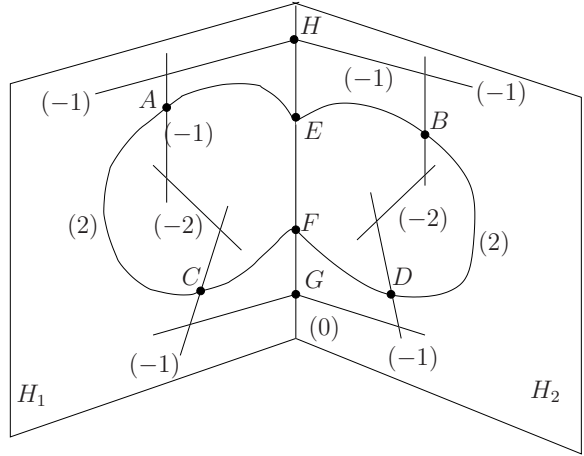,width=3in} 
\end{center}
\caption{The planes $H_{1}$ and $H_{2}$ after flops in the exceptional
  curves on $\widetilde{S}_{3}$.}
\label{fig:H.after.flops} 
\end{figure}

Next we proceed with the flops on page~\pageref{page:flops}.  The
picture for $H_{1}$ and $H_{2}$ after doing these flops but before
performing the small resolutions (at $G$ and $H$) is sketeched in
Figure~\ref{fig:H.after.flops}. 

\

\

\begin{figure}[!ht]
\begin{center}
\psfrag{(-1)}[c][c][0.6][0]{{$(-1)$}}
\psfrag{(0)}[c][c][0.6][0]{{$(0)$}}
\psfrag{(-4)}[c][c][0.6][0]{{$(-4)$}}
\psfrag{(-2)}[c][c][0.6][0]{{$(-2)$}}
\psfrag{(2)}[c][c][0.6][0]{{$(2)$}}
\psfrag{2:1}[c][c][0.6][0]{{2:1}}
\psfrag{A}[c][c][0.6][0]{{$A$}}
\psfrag{B}[c][c][0.6][0]{{$B$}}
\psfrag{C}[c][c][0.6][0]{{$C$}}
\psfrag{D}[c][c][0.6][0]{{$D$}}
\psfrag{E}[c][c][0.6][0]{{$E$}}
\psfrag{F}[c][c][0.6][0]{{$F$}}
\psfrag{G}[c][c][0.6][0]{{$G$}}
\psfrag{S}[c][c][0.6][0]{{$S$}}
\psfrag{S3}[c][c][0.6][0]{{$\widetilde{S}_{3}$}}
\psfrag{P2}[c][c][0.6][0]{{$\mathbb{P}^{2}$}}
\psfrag{F2}[c][c][0.6][0]{{$\mathbb{F}_{2}$}}
\psfrag{A2 singularity}[c][c][0.6][0]{\begin{minipage}[c]{0.8in}{$A_{2}$
      surface \\ singularity}\end{minipage}} 
\epsfig{file=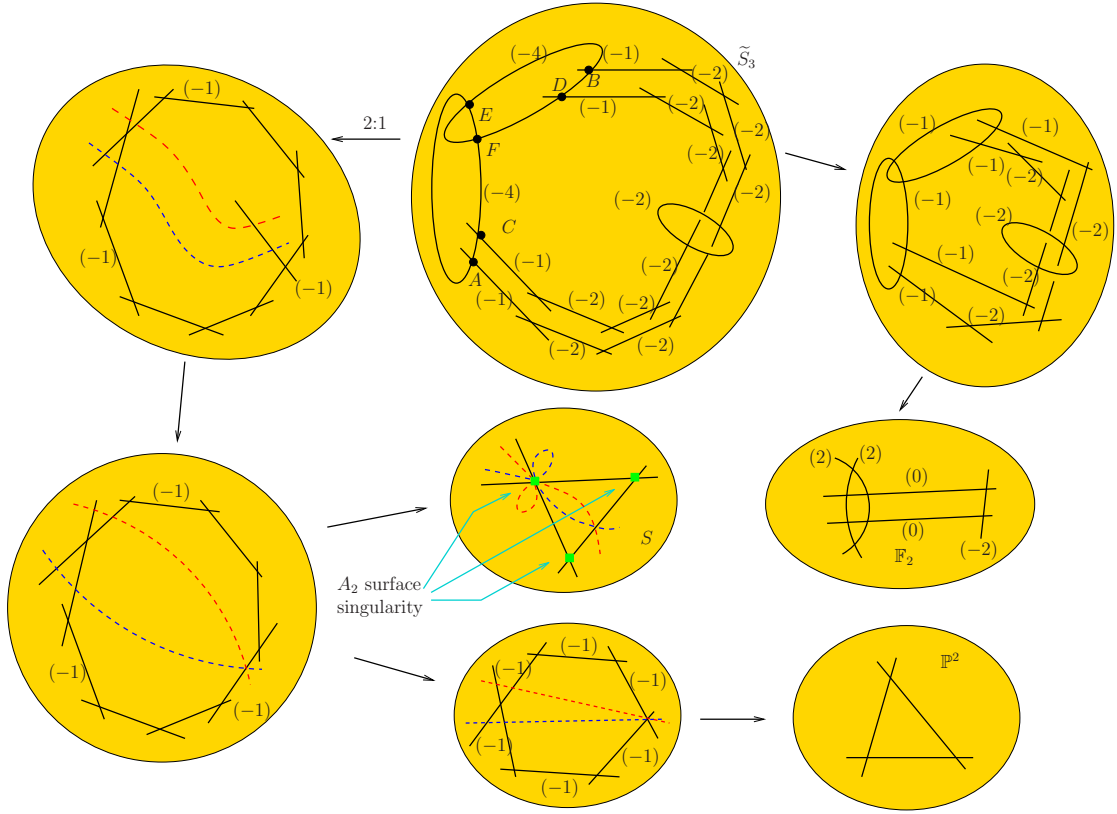,width=5.8in} 
\end{center}
\caption{The effect of the flops on $\widetilde{S}_{3}$. The dashed
  curves depict the branch divisors of the double cover structure on
  $\widetilde{S}_{3}$.}
\label{fig:S3.after.flops} 
\end{figure}

\noindent
It is clear from Figure~\ref{fig:H.after.flops} that the two
components of $Y_{0}$ corresponding to $H_{1}$ and $H_{2}$ can be
obtained by blowing up $\mathbb{F}_2$ at two points of the curve
$q_{3}$. Therefore, these two components look similar to the $F$
component. To get to the final picture in Theorem~\ref{theo-genus2},
we have to perform small resolution at $H$ and $G$ in different
components. In particular the curve $q_{3}$ on the two $\mathbb{F}_2$
will be blown up at 3 points total. To get the picture for $\bfY_{0}$,
we have to blow up points on the conics $q_{1}$ and $q_{2}$ as
well. Specifically we have to blow up $q_{1}$ and $q_{2}$ three times
each at $A$, $C$ and $B$, $D$ respectively. The whole process is
summarized in Figure~\ref{fig:S3.after.flops}. These blow-ups correspond
to flops on the total space $Y$: they can be achieved by flopping
suitable  $(-1)$ curves on the
component $F$. The important fact here is that the
intersection curves $q_1$, $q_2$, and $q_3$, after all resolutions and
flops, meet at $E$ and $F$, and also have two more marked 
points with $(-1)$ curves
at each one on the corresponding surface.

\

\

\ \hfill $\Box$

\

\bigskip

\begin{rem} \label{rem:genus2.Karoubi.revisited} Note that the
description of the fiber $\bfY_{0}$ in the previous theorem shows that
if we vary the symplectic curve $(\underline{C},B+i\omega)$ the
complex algebraic surface $\bfY_{0}$ varies but has one dimensional
moduli parameterized by the cross-ratio parameter $\lambda$, which in
turn depends only on the symplectic volume parameter $a_{2} =
\exp\left( - \int_{\mathbb{P}^{1}\times \op{pt}} (B + i\omega)
\right)$.  This is not surprising since the map from torus equivariant
K\"{a}hler classes on $\mathbb{P}^{1}\times \mathbb{P}^{1}$ has a one
dimensional kernel and so there is a non-trivial relation between
$a_{1}$ and $a_{2}$ when we view them as K\"{a}hler parameters on $C$.

The formal neighborhood of the singular locus of $\bfY_{0}$ likely
depends of the gluing parameter $\lambda$ as well. Through
the localization theorem \cite[Theorem~2.10]{orlov-completion} this
gives and apparent dependence on $\lambda$ for the
Karoubi closure
\[
\overline{D_{\op{sing}}}(\bfY_{0}) =
\overline{D^{b}}(\bfY_{D},\bw).
\] 
This observation suggests that the conjectural ($A$ to $B$)
equivalence
\[ \label{eq:AtoB} \tag{$A$ to $B$}
\overline{\Fuk}(\underline{C},B+i\omega) \cong
  \overline{D^{b}}(\bfY_{D},\bw)
\]
predicted by HMS is really a family of different equivalences labeled
  by $\lambda$. 

One can try to detect the $\lambda$-dependence in \eqref{eq:AtoB} 
on the level of $K$-groups. From the work of Abouzaid
\cite{abouzaid} it is known that the $K$-group of the Karoubi-closed
derived Fukaya category $\overline{\Fuk}(\underline{C},B+i\omega)$ is
not finitely generated. More precisely from \cite{abouzaid} it follows
that the $K$-group of the Fukaya category is given by
\[
\begin{split}
K_{0}(\overline{\Fuk}(\underline{C},B+i\omega)) & =
H_{1}(\underline{C},\mathbb{Z})\oplus \mathbb{Z}/\chi(\underline{C})
\oplus \mathbb{C}^{\times} \\
& = \mathbb{Z}^{4}\oplus \mathbb{Z}/2\oplus \mathbb{C}^{\times}.
\end{split}
\]
Combining this with the recent proof by Seidel \cite{SEID} of the
mirror equivalence \eqref{eq:AtoB} we conclude
that
\[
K_{0}(\overline{D^{b}}(\bfY_{D},\bw)) \cong \mathbb{Z}^{4}\oplus
\mathbb{Z}/2\oplus \mathbb{C}^{\times},
\]
but perhaps there is a whole family of such isomorphisms which depends
non-trivially on $\lambda$.

In fact in the last formula can be verified by a direct calculation
\cite{dima-high} in $\overline{D^{b}}(\bfY_{D},\bw)$ but we will skip
this verification here. Some evidence for that is provided by the fact
that the $K$-group of the uncompleted category $D^{b}(\bfY_{D},\bw) =
D^{b}_{\op{sing}}(\bfY_{0})$ depends \cite{dima-high} non-trivially on
the moduli of $\bfY_{0}$. 

Secifically suppose that $\lambda$ is a primitive
$n$-th root of unity and that $a_{2} \in \mathbb{C}$ is such that
$\lambda$ is the ratio of the roots $x^{2} + x + a_{2} = 0$ in some
order. (In other words we have $a_{2} = \lambda/(\lambda + 1)^{2}$.)
If $\bfY_{0}$ is the singular surface of modulus $\lambda$
constructed in theorem~\ref{theo-genus2-symmetric}, then one can check
\cite{dima-high} that
\[
K_{0}(D^{b}_{\op{sing}}(\bfY_{0})) \cong \mathbb{Z}^{4}\oplus \mathbb{Z}/2n.
\]
\end{rem}

\subsection{General hyperelliptic curve}

Let us consider a hyperelliptic curve $C$ of genus $k-1$. We can
think of it as a difisor of bi-degree $(2, k)$ in the Hirzebruch
surface ${\mathbb F}_0={\mathbb P}({\mathcal O}_{{\mathbb
P}^{1}}\oplus {\mathcal O}_{{\mathbb P}^{1}})$ . Similarly to the
the genus two case $k = 3$ we get the weights of the torus action
in the table:

\

\begin{center}
\begin{tabular}{|c|c|c|c|c|    } \hline
$z_{1}$ & $z_{2}$ &                     $w_{1}$ & $w_{2}$ & $p$ \\
\hline \hline
$1$ & $1$ &             $0$ & $0$ & $- k$ \\ \hline
$0$ & $0$ &             $1$ & $1$ & $- 2$ \\ \hline

\end{tabular}
\end{center}

\

The Hori-Vafa procedure \cite{HV} gives an affine mirror
$(Y^{\op{aff}},\bw^{\op{aff}})$ of $(\underline{C},B+i\omega)$ where 
$Y^{\op{aff}}$ is the threefold in $\mathbb{C}^{\*5}$ defined by the
equations
\begin{equation}
\label{eq:affineY.genus.k-1}
Y^{\op{aff}} \ : \ 
\left\{ \begin{split}
x_1\cdot x_2  = a_1\cdot x_3^{k} \\ 
x_4\cdot x_5  = a_2\cdot x_3^{2}.
\end{split} \right.
\end{equation}
and the pencil  $\bw^{\op{aff}} : Y^{\op{aff}} \to \mathbb{C}$ is
defined by the formula 
\begin{equation} \label{eq:affinew.genus.k-1}
\bw^{\op{aff}} = x_1 + \cdots + x_5.
\end{equation}
The corresponding partial compactification and resolution can be
studied in a way similar to the genus two case. We omit the detailed
description of the singular fibers of the resulting $\bw : Y \to
\mathbb{C}$. For completeness we give the combinatorial description of
the critical locus of $\bw$.

\begin{figure}[!ht]
\begin{center}
\epsfig{file=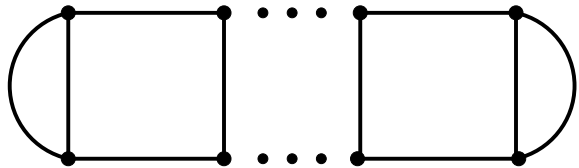,width=3in} 
\end{center}
\caption{Hyperelliptic curves}
\label{fig-wk} 
\end{figure}

\begin{theo} \label{theo:shape.genus.k-1} Let $(Y,\bw)$ be a
  Landau-Ginzburg mirror of 
  $(\underline{C},B+i\omega)$ which compactifies and resolves the
  affine mirror given by \eqref{eq:affineY.genus.k-1} and
  \eqref{eq:affinew.genus.k-1}. Then:
\begin{itemize}
\item $\bw$ has exactly $k$ critical  values;
\item The zero dimensional part of the critical set of $\bw$ consists of
$k-1$ isolated points of $Y$, the corresponding fibers of $\bw$ have
  simple $A_{1}$ singularities at those points;  
\item The positive-dimensional part of the critical locus of $\bw$
  sits over the critical value $0$ and consists of $3k-6$ rational
  curves forming a degenerate curve of genus $k-1$ whose dual graph
  is shown in figure \ref{fig-wk}.
\end{itemize}
\end{theo} 

\

\

\noindent
Recently Seidel \cite{SEID} and \cite{efimov} Efimov proposed an
alternative construction of the Hori-Vafa mirror of the hyper-elliptic
curve $C$. The idea is not to treat $C$ as a curve emebedded in a
toric surface but rather as an etale $\mathbb{Z}/(2k-1)$-Galois cover
of a genus zero orbifold. Understanding the mirror of such an orbifold
and tracing the effect of the Galois action leads to an affine
Landau-Ginzburg model whos total space has quotient
singularities. These in turn can be resolved by introducing
appropriate orbifold structures - a procedure which is much more
tracktable than the toric resolution procedure we employed in the
previous section. Not surprisingly the critical locus of the
superpotential on the orbifold Landau-Ginzburg mirror has the same
combinatorial structure as the one we obtained in
Theorem~\ref{theo:shape.genus.k-1}.

By analyzing this geometry and combining it with a version of the
McKay correspondence with potentials Efimov proves \cite{efimov} the
($A$ to $B$) equivalence
\[ 
 \tag{$A$ to $B$}
\overline{\Fuk}(\underline{C},B+i\omega) \cong
  \overline{D^{b}}(Y_{D},\bw)
\]
predicted by HMS. Here as before $D$ is a sufficiently small disk
in $\mathbb{C}$ which only contains the critical value $0$ of $\bw$.

\

\noindent
The above orbifold procedure is not directly compatible with our
construction of the mirror. However, as it was pointed out to us by
P.Seidel, the trick with orbifold resolutions can be easily modified
to fit with our birational model. We will sketch this modification
here since the resulting geometry ties very nicely with the approach
of Seidel and Efimov. Start again with the affine Hori-Vafa mirror
$(Y^{\op{aff}},\bw^{\op{aff}})$ given by \eqref{eq:affineY.genus.k-1}
and \eqref{eq:affinew.genus.k-1}. Before we compactify the fibers of
$\bw^{\op{aff}}$ and start resolving the singularities of the
resulting total space we can simplify the problem by first resolving
the singularities of $Y^{\op{aff}}$ and then choosing a partial
compactification compatible with this resolution. This can be done in
two steps - first pass to an orbifold cover of $Y^{\op{aff}}$ with
mild singularities, and then construct a crepant resolution of that
cover. Concretely, consider the quadric cone $Z^{\op{aff}} \, : \,
z_{1}z_{2} = z_{4}z_{5}$ sitting in the $\mathbb{C}^{4}$ with
coordinates $z_{1}$, $z_{2}$, $z_{4}$, $z_{5}$.  Consider also the
natural map $Z^{\op{aff}} \to Y^{\op{aff}}$ given by
$(z_{1},z_{2},z_{4},z_{5}) \mapsto
(z_{1}^{k},z_{2}^{k},z_{1}z_{2},z_{4}^{2},z_{5}^{2})$. The map
$Z^{\op{aff}} \to Y^{\op{aff}}$ is a ramified $(\mathbb{Z}/k)\times
(\mathbb{Z}/2)$-Galois cover, where $(\lambda,\mu) \in
(\mathbb{Z}/k)\times (\mathbb{Z}/2)$ acts on $Z$ by
$(\lambda,\mu)\cdot (z_{1},z_{2},z_{4},z_{5}) := (\lambda z_{1},
\lambda^{-1} z_{2}, \mu z_{4}, \mu^{-1}z_{5})$.  Choose now a small
resolution $\widehat{Z} \to Z^{\op{aff}}$ of this quadratic
cone which is compatible with the group action.  For instance we can
take $\widehat{Z}$ to be the blow-up of $Z$ along the plane $z_{1} =
z_{4} = 0$. Note that by construction $\widehat{Z}$ embeds in
$\mathbb{C}^{5}\times \mathbb{P}^{1}$ and this embedding is
equivariant for the natural $(\mathbb{Z}/k)\times
(\mathbb{Z}/2)$-actions. Proceeding as before we can partially
compactify $\widehat{Z}$ by closing it inside $\mathbb{P}^{4}\times
\mathbb{C}\times \mathbb{P}^{1}$. The resulting space $Z$ is smooth
and the quotient variety $Z/((\mathbb{Z}/k)\times (\mathbb{Z}/2))$ is
a partial compactification and a small modification of
$Z^{\op{aff}}$. By Hartogs' theorem the potential $\bw^{\op{aff}}$
induces a proper holomorphic function $\bff : Z\to \mathbb{C}$ and the
group $(\mathbb{Z}/k)\times (\mathbb{Z}/2)$ acts along the fibers of
$\bff$. In particular $\bff$ will descend to a function (which we
denote again by $\bff$) on the stack quotient $\mathcal{Z} :=
\left[Z/((\mathbb{Z}/k)\times (\mathbb{Z}/2))\right]$ and so the
orbifold Landau-Ginzburg model $(\mathcal{Z},\bff)$ can be viewed as a
stacky crepant resolution of $(Y^{\op{prop}},\bw^{\op{prop}})$. The
categorical McKay correspondence of \cite{BKR} and \cite{kawamata}
applies to this situation and allows us to identify the derived
categories of $\mathcal{Z}$ and the crepant resolution $Y$.  Combined
with our argument from section~\ref{sec:flops} this implies that we
have an equivalence of categories $D^{b}(Y,\bw) \cong
D^{b}(\mathcal{Z},\bff)$.

This interpretation of the $B$-model category has some definite
computational advantages. First, since $\mathcal{Z}$ is a quotient
stack, we can compute $D^{b}(\mathcal{Z},\bff)$ as the
$(\mathbb{Z}/k)\times (\mathbb{Z}/2)$-equivariant derived category of
the potential $\bff : Z \to \mathbb{C}$. Also, since the critical
locus of $\bff$ is disjoint from the boundary $Z-\widehat{Z}$ we can
use the localization theorem from \cite{O} to argue that
$D^{b}(\mathcal{Z},\bff)$ is equivalent to the $(\mathbb{Z}/k)\times
(\mathbb{Z}/2)$-equivariant derived category of the potential $\bff :
\widehat{Z} \to \mathbb{C}$. The later category should be easy to
compute since $\widehat{Z}$ is just the total space of
$\mathcal{O}(-1)\oplus \mathcal{O}(-1)$ on $\mathbb{P}^{1}$ and so its
derived category is computed by a quiver. The superpotential $\bff$
gives an $A_{\infty}$-deformation of this algebra and smashing this
deformation with group $(\mathbb{Z}/k)\times
(\mathbb{Z}/2)$ we get an explicit $A_{\infty}$ algebra that computes
the $B$-model category $D^{b}(\mathcal{Z},\bff)$.

We also have a similar orbifold picture for the symplectic
hyperellyptic curve $(\underline{C},B+i\omega)$ on the
$A$-side. Indeed, we can view $\underline{C}$ as an etale
$(\mathbb{Z}/k)\times (\mathbb{Z}/2)$-Galois cover of an orbifold
$S^{2}$ with four orbifold points with stabilizers $\mathbb{Z}/k$,
$\mathbb{Z}/k$, $\mathbb{Z}/2$, and $\mathbb{Z}/2$ respectively. Using
the techniques of Seidel and Efimov we can describe the Fukaya
category of such a $2k$-sheeted  cover in terms of curves on the four
punctured sphere. By repeating the analysis in \cite{efimov}
one should be able to match this description of the Fukaya category
directly with the $B$-side deformed smash algebra that we  constructed
above. Note that the orbifold $\mathbb{P}^{1}$ with four
orbifold points of indices $(k,k,2,2)$ appears naturally on the
$B$-side as well. Specifically the orbifold
$\left[\widehat{Z}/((\mathbb{Z}/k)\times (\mathbb{Z}/2))\right]$ can
be identified with the total space of a vector bundle $L\oplus L$ on
this orbifold $\mathbb{P}^{1}$, where $L$ is the unique theta
characteristic on the orbifold $\mathbb{P}^{1}$. 

The above orbifold picture provides yet another description of the
mirror of the hyperelliptic curve and gives an alternative setting for
proving Efimov's theorem.

\begin{rem} The shape of the critical locus of $\bw$ suggests that both
  $D^{b}(Y_{D},\bw)$ and \linebreak $\Fuk(\underline{C},B+i\omega)$ may be
  amenable to cutting and pasting computations. In particular it will
  be interesting to try and compute $\Fuk(\underline{C},B+i\omega)$ by
  gluing $\underline{C}$ small pieces - e.g. by using a pair of pants
  decomposition. It will also be interesting to see how the mapping class
  group interacts with the mirror functor.
\end{rem}

\

\noindent
We will leave this question for the future. For now we proceed with
the analysis of the mirror of the intersection of two quadrics which is
the other half of our conjecture.

\section{Intersection of two quadrics in $\CP^5$} \label{sec-quadrics}

In this section, we give a detailed description of the Landau-Gizburg
mirror $\bw : Y \to \mathbb{C}$ of the Fano three-fold $X_{4}$ of
degree 4 in $\CP^5$. We will show that the topological structure of
the critical set of $\bw$ is identical with the one corresponding to
the genus two curve $C$. Assumimg the HMS conjecture, we also give
evidence supporting our Conjecture~\ref{con:side.A.equiv}, which in
this case predicts that under an appropriate identification of
complexified K\"{a}hler classes the Fukaya category of $\underline{C}$
maps fully faithfully to the Fukaya category of the three-fold
$\underline{X}_{4}$.

The Fano manifold $X_4$ is a complete intersection of two quadrics in
$\CP^5$. If we equip the underlying $C^{\infty}$-manifold
$\underline{X}_{4}$ with the unique up-to-scale complexified
K\"{a}hler form which is a restriction of a torus uinvariant form on
$\mathbb{P}^{5}$, then according to \cite{HV} its affine
Landau-Ginzburg mirror is given by the subvariety in $\mathbb{C}^{6*}$
cut out by the system of equations:
\[
Y^{\op{aff}} \ : \
\left\{
\begin{split}
x_1x_2x_3x_4x_5x_6  & = 1 \\
x_1 + x_2 & = -1  \\
x_3 + x_4 & = -1  
\end{split} \right.
\]
and equippped with the superpotential $\bw^{\op{aff}} = x_5 + x_6$.

To construct the partial coompactification $Y^{\op{prop}}$ 
of $Y^{\op{aff}}$ we
consider the
closure of this variety in \,\,$\CP^1_{(x_0:x_1)}\times
\CP^1_{(y_0:y_1)}\times \CP^1_{(z_0:z_1)}\times \mathbb{C}_w$\,\,
given by the equation
\[
Y^{\op{prop}} \ : \ 
 x_1(x_0 - x_1)y_1(y_0 - y_1)z_1(wz_0 - z_1) = x_0^2y_0^2z_0^2.
\]
Here the coordinates of each $\mathbb{P}^{1}$-factor are indicated in the
subscripts. As before the extension $\bw^{\op{prop}}$ of the potential
$\bw^{\op{aff}}$ is given by projecting $Y^{\op{prop}}$ onto the $w$-line.

The desingurization $Y$ of this variety and the structure of the
singular fibers of the resulting potential $\bw$ are summarized in the
following theorem.

\begin{theo} \label{theo-quadrics} There exists a crepant resolution 
$Y$ of $Y^{\op{prop}}$ so that the fibers of the function $\bw : Y \to
  \mathbb{C}$ induced from $\bw^{\op{prop}}$ can be described as follows:
\begin{itemize}
\item For $w\neq -8, 0, 8 $, the fiber $Y_w := \bw^{-1}(w)$ is
smooth and has a structure of an elliptic fibration over $\CP^1$ with
five singular fibers: one of type $I_8$, two of type $I_1^*$, and two
of type $I_1$. Hence, $Y_w$ is an elliptic $K3$ surface with Picard
rank 19. 
\item For $w = \pm 8$ the fiber $Y_{w}$ is singular with an isolated $A_{1}$
  singularity. Its minimal resolution is an elliptic $K3$ with an
  $I_{8} + 2\cdot I_{1}^{*} + I_{2}$ configuration of singular fibers.
\item The central fiber $Y_0$ of the
fibration $w : Y \to \CP^1$ is a union of three rational surfaces
$Y = \widetilde{S}_1 \cup \widetilde{S}_2 \cup
\widetilde{S}_{3}$. These three surfaces intersect in three 
copies of $\mathbb{P}^{1}$: 
\[
q_1 = \widetilde{S}_1 \cap \widetilde{S}_{3}, \quad q_2 =
\widetilde{S}_2 \cap \widetilde{S}_{3}, \quad \text{and} \quad q_3 =
\widetilde{S}_1 \cap \widetilde{S}_2.
\] 
The curves $q_1$, $q_2$, and $q_3$ intersect in two points,
$R_1$ and $R_2$. Hence, $Y_0$ is a type III deformation of a K3
surface. The geometry of the surfaces and the position of the curves
on them is explained as follows.
\begin{enumerate}
\item[(i)] To construct $\widetilde{S}_{3}$, choose four generic
  rulings on $\mathbb{F}_0$ which intersect at points $P_1$, $P_2$,
  $Q_1$, and $Q_2$ so that $P_1$ and $Q_1$ are "diagonal" points
  w.r.t. the chosen rulings. Choose further two generic $(1, 1)$
  divisors, $d_1$ and $d_2$ such that $d_1$ passes through $P_1$ and
  $Q_1$, and $d_2$ passes through $P_2$ and $Q_2$. Let $d_1 \cap d_2 =
  \{R_1, R_2\}$. The component $\widetilde{S}_{3}$ is the blow up of
  $\mathbb{F}_0$ in points of depth two on $d_1 \cup d_2$ and support at
  $P_1, P_2, Q_1, Q_2$. The proper transforms of the divisors $d_1$
  and $d_2$ are the curves $q_1$ and $q_2$ on $\widetilde{S}_{3}$.

\item[(ii)] The set up for constructing $\widetilde{S}_1\cup
  \widetilde{S}_2$ is the same as in (i). Only now we blow up
  $\mathbb{F}_0$ in the vertices of the quadrilateral only once;
  consider two copies of this modified surface and identify them along
  the proper transform of $d_1$ so that the points $R_1$ and $R_2$
  from the two copies of the surface match; then we blow up either the
  first or the second component at two points on $d_1$ but different
  from $R_1$ and $R_2$. The resulting variety is $\widetilde{S}_1 \cup
  \widetilde{S}_2$. The two divisors $d_2$ correspond to $q_1$ and
  $q_2$ on $\widetilde{S}_1$ and $\widetilde{S}_2$ respectively. The
  proper transform of $d_1$, on either of $\widetilde{S}_i$, is the
  curve $q_3$ on the corresponding component of $Y_0$.

This way of constructing $\widetilde{S}_1 \cup \widetilde{S}_2$ is
clarified by the fact that, before we blow up the last two points on
$d_1$, there is an obvious $\ZZ_2$ symmetry of the union of the three
components which identifies the two exceptional ones as we do
here. According to this symmetry, for the case of the third component,
we have to choose on the corresponding $\mathbb{F}_0$ the configuration
of diagonal divisors symmetric with respect to the $\ZZ_2$ symmetry
which exchanges the members of one of the rulings.

\end{enumerate}
The rank of the Picard group of $Y_0$ is 21.
\end{itemize}
\end{theo}

\noindent
{\bfseries Proof:} Consider the intersection $X_4 = Q_1 \cap Q_2
\subset \CP^5$ of two generic quadrics in $\CP^5$. The space $X_{4}$
has a unique up to scale complexified symplectic form which is a
restriction of a torus invariant complexified symplectic form on
$\mathbb{P}^{5}$. For concretensess we normalize this form to be $2\pi
i \omega_{FS}$ where $\omega_{FS}$ is the Fubini-Studi form on
$\mathbb{P}^{5}$. 

The Hori-Vafa recipe \cite{HV} identifies the affine Landau-Ginzburg
mirror of $(\underline{X}_{\, 4}, 2\pi i \omega_{FS})$ with the family
$\bw^{\op{aff}} : Y^{\op{aff}} \to \mathbb{C}$ where
\[
Y^{\op{aff}} \ : \
\left\{
\begin{split}
x_1x_2x_3x_4x_5x_6  & = 1 \\
x_1 + x_2 & = -1  \\
x_3 + x_4 & = -1  
\end{split} \right.
\]
and the superpotential is given by $\bw^{\op{aff}} = x_5 + x_6$.
Note that any rescaling of the Fubini-Studi form will only affect the
affine mirror by introducing some non-zero complex constant $r$ instead of
the $1$ in the right-hand-side of the first equation. The resulting
affine varieties with potentials are clearly isomorphic (say by
scaling $x_{5}$ and $x_{6}$ by the square root of $r$). Thus the
scaling parameter for the symplectic volume of $X_{4}$ is not a
modulus of the mirror Landau-Ginzburg model and we can work with the
normalized mirror above without any loss if generality.

We consider the closure of this variety in
$\CP^1_{(x_0:x_1)}\times \CP^1_{(y_0:y_1)}\times
\CP^1_{(z_0:z_1)}\times \C_w$ given by the equation
\[
Y^{\op{prop}} \ : \  x_1(x_0 - x_1)y_1(y_0 - y_1)z_1(wz_0 - z_1) =
x_0^2y_0^2z_0^2 
\]
A crepant resolution $Y$ of this variety together with the induced
superpotential will be Landau-Ginzburg mirror fibration of the Fano
manifold $(\underline{X}_{\, 4}, 2\pi
i \omega_{FS})$. The
singular locus of $Y^{\op{prop}}$ consists of the following 
one dimensional components: 
\[
\begin{array}{lll}
l_1: [(0:1);(1:0);(1:0);w] \quad & m_1: [(1:0);(0:1);(1:0);w] \quad \\ 
l_2: [(0:1);(1:0);(1:w);w] \quad & m_2: [(1:0);(0:1);(1:w);w] \quad \\
l_3: [(0:1);(1:1);(1:0);w] \quad & m_3: [(1:1);(0:1);(1:0);w] \quad \\
l_4: [(0:1);(1:1);(1:w);w] \quad & m_4: [(1:1);(0:1);(1:w);w] \quad \\
n_1: [(1:0);(1:0);(0:1);w] \quad &  p_1: [(0:1);(y_0:y_1);(1:0);0]
\quad \\
n_2: [(1:0);(1:1);(0:1);w] \quad & p_2: [(x_0:x_1);(0:1);(1:0);0]
\quad \\
n_3: [(1:1);(1:0);(0:1);w] \quad & \\
n_4: [(1:1);(1:1);(0:1);w] \quad & 
\end{array}
\]
All these components have transversal singularity of type $A_1$, and,
at generic points, can be resolved by blowing up the components in the
ambient space. The fiber $Y^{\op{prop}}_w$, for $w \neq 0$, intersects
transversally 12 of the lines above. The fiber $Y_w$ of the
desingularization $Y$ is then the resolution of the $12(A_1)$
singularities of $Y^{\op{prop}}_w$.

Again we begin with a description of the general fiber.

\begin{claim}
For $w\neq -8, 0, 8 $, $Y_w$ is a smooth $K3$ surface and has a
structure of an elliptic fibration over $\CP^1$ with five singular
fibers: one of type $I_8$, two of type $I_1^*$, and two of type
$I_1$. The rank of the Picard group of $Y_w$ is 19. The fibers $Y_{\pm
8}$ are singular with an isolated $A_{1}$ singularity. The minimal
resolution of either of these two fibers is an elliptic $K3$ with an
$I_{8} + 2\cdot I_{1}^{*} + I_{2}$ configuration of singular fibers.
\end{claim}
{\bf Proof.} \ Observe that $Y^{\op{prop}}$ is a double cover of
$S'\times \C_w $ where
\[
S': \{ x_1(x_0 - x_1)y_1(y_0 - y_1)T_1 = x_0^2\,y_0^2\,T_0\} \subset
\CP^1_{(x_0:x_1)}\times \CP^1_{(y_0:y_1)}\times \CP^1_{(T_0:T_1)}
\]
and the cover map $\varphi': Y^{\op{prop}} \to S'\times \C_w$ is
induced by the double cover $\psi: \CP^1_{(z_0:z_1)}\to
\CP^1_{(T_0:T_1)}$ where $\psi ((z_0:z_1)) = (z_0^2: z_1(w\,z_0 -
z_1))$.  The surface $S'$ has four singular points of type $A_1$ and
with respect to the projection to $\CP^1_{(T_0:T_1)}$, is a rational
elliptic fibration with singular fibers of type $I_4$, at $T_0 = 0$,
type $I_1$ at $T_1 = 16T_0$, and a double fiber of nodal type, at
$T_1= 0$. The four $A_1$ singularities of $S'$ belong to this latter
fiber (neither of these is a nodal point there).  The restriction of
$\varphi'$ to $Y^{\op{prop}}_w$ is a double cover map $\varphi'_w:
Y^{\op{prop}}_w \to S'$ The branching for $\psi$ occurs at $T_0 = 0$
and $T_1 = (w/2)^2\,T_0$. Hence, the branch divisor of $\varphi'_w$
consists of two fibers of $S'$: at $T_0 = 0$ and at $T_1 =
(w/2)^2\,T_0$. This implies that, for $w \neq -8, 0, 8$,\,\,
$Y^{\op{prop}}_w$ has a structure of an elliptic fibration over
$\CP^1_{(T_0:T_1)}$ with singularities induced by those of $S'$, 8
total, and of the self-intersection of the branch divisor, 4
total. These are resolved by blowing up the lines $\,\,l_i, m_i, n_i,
\,(i = 1,\dots ,4)\,\,$, and the result is an elliptic fibration $Y_w
\to \CP^1$ with an $I_{8} + 2\cdot I_{1}^{*} + 2I_{1}$ configuration of
singular fibers.

Denote by $S$ the resolved $S'$ blown up additionally in the nodes of
its $I_4$ fiber. Then $S$ is a rational elliptic fibration with three
singular fibers. One of them looks like an $I_8$ fiber but has has
four disjoint components of self-intersection $-4$ and other four of
self-intersection $-1$. The second singular fiber is of type
$I_1^*$. The third is a type $I_1$ fiber. The surface $Y_w$ is a
double cover of $S$ branched in a union of a smooth fiber and of the
four $-4$ lines of the $I_8$-like fiber It is easy to see now that
there are only 24 $(-2)$-lines on $Y_w$. The fibration structure on
$Y_w$ gives two relations between these $(-2)$-lines in
$\op{Pic}(Y_w)$. If we repeat, word for word, our considerations for
the projections of $Y^{\op{prop}}_w$ to the rest of the $\CP^1$s (with
coordinates $(x_0:x_1)$ and $(y_0:y_1)$), we will get two more ways to
express $Y_w$ as an elliptic fibration with the same types and numbers
of singular fibers. This way we get 4 more relations between the
$-2$-lines in $\op{Pic}(Y_w)$. A direct check shows that among the six
relations we get this way only five are linearly independent. Since
the Picard group on such a K3 surface is generated by the
$(-2)$-components of the fibers, we get that the rank of this group
for $Y_w$ is no more than 19. By a direct chase one can find 19
linearly independent such lines in our case.

The fibers $Y_{\pm 8}$ are branched over $S$ at the four $-4$-lines as
before, and at the $I_1$ fiber of $S$. This way we immediately get the
structure of $Y_{\pm 8}$.  This completes the proof of the claim. \
\hfill $\Box$ 

\

\bigskip

\noindent
Next we turn to the central fiber $Y^{\op{prop}}_0$ of
$Y^{\op{prop}}$. This fiber will be affected also by the resolution of
the lines $p_1$ and $p_2$, which are contained in this fiber. Along
these two lines, $Y^{\op{prop}}_0$ has a self-intersection of nodal
type and resolving them will produce two extra components in
$Y_{0}$. The double cover $Y^{\op{prop}} \to S'$ will help us
understand the structure of $Y_{0}$. This cover is branched along the
fiber $I_4$ as well as in the double nodal fiber (for $T_1 = 0$). The
double portion of the branch divisor is responsible for the non-normal
singularities of $Y^{\op{prop}}_0$, and to resolve these is equivalent to
(partially) "deleting" the doble-lines part from the branch divisor.

To properly carry out this "deletion", we start by blowing up the line
$p_1$ given by $x_0 = z_1 = w = 0$. By setting
\[
(w: x_0/x_1: z_1/z_0) = (u_0: u_1: u_2)
\]
we get the blown up $Y^{\op{prop}}$ in the form
\[
(x_0/x_1 -1)\cdot y_1\cdot(y_0 - y_1)\cdot u_2\cdot(u_0 - u_2) =
u_1^2\cdot y_0^2 
\]
where $x_0/x_1 = 0$ gives the exceptional divisor, and the proper
transform of $Y^{\op{prop}}_0$ is given by $u_0 = 0$.

\begin{figure}[!ht]
\begin{center}
\psfrag{F_0}[c][c][0.8][0]{{$\mathbb{F}_{0}$}}
\psfrag{u_1=0}[c][c][0.8][0]{{$u_{1}=0$}}
\psfrag{u_2=0}[c][c][0.8][0]{{$u_{2}=0$}}
\psfrag{u_0 - u_2 = 0}[c][c][0.8][0]{{$u_{0}-u_{2}=0$}}
\psfrag{P^2}[c][c][0.8][0]{{$\mathbb{P}^{2}$}}
\psfrag{(-1)}[c][c][0.8][0]{{$(-1)$}}
\psfrag{(2)}[c][c][0.8][0]{{$(2)$}}
\psfrag{(0)}[c][c][0.8][0]{{$(0)$}}
\psfrag{(-2)}[c][c][0.8][0]{{$(-2)$}}
\psfrag{Branch locus}[c][c][0.8][0]{{Branch locus}}
\psfrag{(lambda:mu)}[c][c][0.8][0]{{$(\lambda:\mu)$}}
\psfrag{(1:4a_1)}[c][c][0.8][0]{{$(1:4a_{1})$}}
\psfrag{mu = 0}[c][c][0.8][0]{{$\mu=0$}}
\psfrag{A_1}[c][c][0.8][0]{{$A_{1}$}}
\psfrag{D_1}[c][c][0.8][0]{{$D_{1}$}}
\psfrag{q_1}[c][c][0.8][0]{{$q_{1}$}}
\psfrag{q_3}[c][c][0.8][0]{{$q_{3}$}}
\psfrag{y_1=0}[c][c][0.8][0]{{$y_{1}=0$}}
\psfrag{y_0=0}[c][c][0.8][0]{{$y_{0}=0$}}
\psfrag{y_0=y_1}[c][c][0.8][0]{{$y_{0}=y_{1}$}}
\psfrag{2:1}[c][c][0.8][0]{{$2:1$}}
\epsfig{file=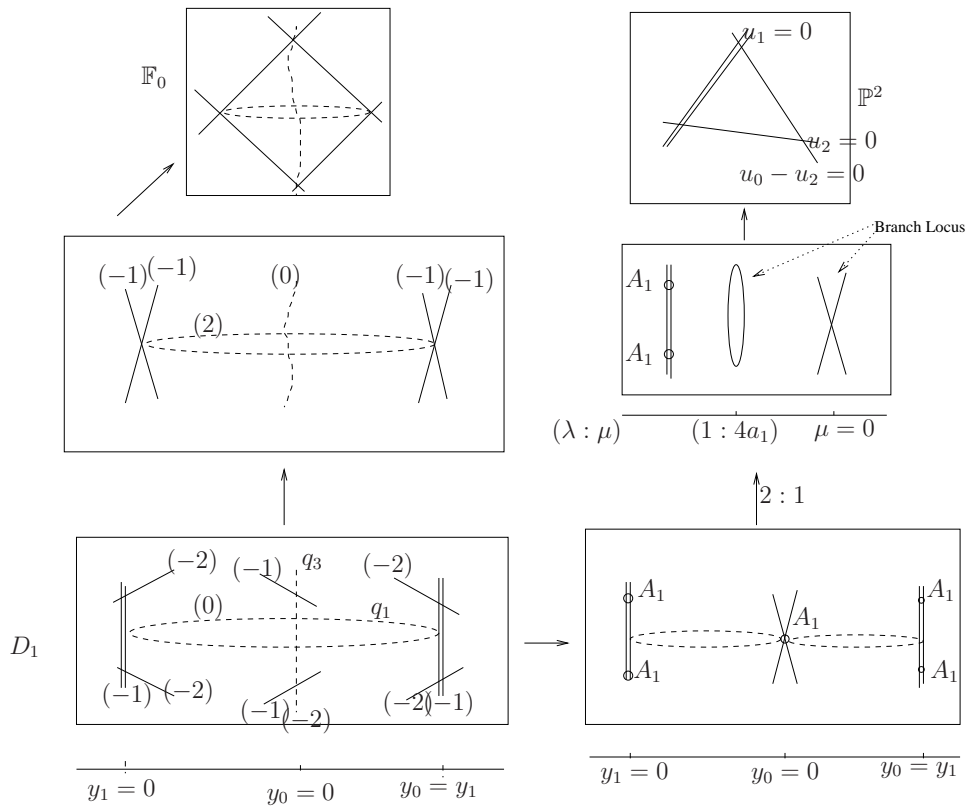,width=5in} 
\end{center}
\caption{The desingularization of $D_{1}$.}
\label{fig:D1} 
\end{figure}

It is easy to check (see Figure~\ref{fig:D1}) that the exceptional
divisor, $D_1$, has five singular points of type $A_1$:
\[
\begin{split}
u_1 & = y_1 = u_2 = 0; \\ 
u_1 & = y_1 = u_0 - u_2 = 0; \\  
u_1 & = y_0 - y_1 = u_2 = 0; \\ 
u_1 & = y_0 - y_1 = u_0 - u_1 = 0; \\ 
y_0 & = u_2 = u_0 = 0.
\end{split}
\] 
The first four singularities get resolved by blowing up the lines
$l_i, \,\,i = 1,\cdots, 4$ while the last one gets resolved by blowing
up the line $p_2$. 

\

\

Before desingularizing $D_1$, we can realise it as a double cover as
follows. Consider the plane $\CP^2$ with coordinates $(u_0: u_1:
u_2)$, and the pencil of conics in it given by the equation
$\lambda\cdot u_2(u_0 - u_1) = \mu\cdot u_0^2$. This pencil defines a
fibration of conics over $\CP^1$ (the line has coordinates $(\lambda:
\mu)$). Then the equations $\lambda = y_1(y_0 - y_1)$ and $ \mu =
y_0^2$ define $D_1$ as a double cover over that fibration. It is
straightforward to see that the desingularization of $D_1$ together
with the curves of intersection with the rest of the components of the
central fiber of $Y$ can be obtained from a Hirzebruch surface
$\mathbb{F}_0$ by blowing it up eight times in the following
centers. Choose a quadrilateral of rulings on $\mathbb{F}_0$ and two $(1,
1)$- divisors which intersect the quadrilateral diagonally. Blow up
$\mathbb{F}_2$ in the vertices of that quadrilateral first, and then in
the two points where the proper transform of one of the diaginal
divisors intersect two of the exceptional curves. This is how the
resolved $D_1$ looks like. The proper transform of the diagonal blown
up twice, $q_1$, is the intersection of $D_1$ with $F$ while the
proper transform of the diagonal blown up four times, i.e. $q_3\,\,$,
is the curve of intersection of $D_1$ with $D_2$.

The second exceptional component of $Y_0$, due to the blow up in
$p_2\,$, can be treated the same way (here we have to blow up the
lines $m_i,\,\, i = 1,\cdots, 4\,$). The only difference is that here
we have to blow up only once the vertices of the quadrilateral. To make
the picture more symmetric, we can flop $D_1$ in a $-1$-curve
intersecting $q_3$.  

\begin{figure}[!ht]
\begin{center}
\psfrag{F_0}[c][c][0.8][0]{{$\mathbb{F}_{0}$}}
\psfrag{(-1)}[c][c][0.8][0]{{$(-1)$}}
\psfrag{(-2)}[c][c][0.8][0]{{$(-2)$}}
\psfrag{(-4)}[c][c][0.8][0]{{$(-4)$}}
\psfrag{Branch pts}[c][c][0.8][0]{{Branch points}}
\psfrag{q_1}[c][c][0.8][0]{{$q_{1}$}}
\psfrag{q_2}[c][c][0.8][0]{{$q_{2}$}}
\psfrag{S}[c][c][0.8][0]{{$S$}}
\psfrag{S'}[c][c][0.8][0]{{$S'$}}
\psfrag{F}[c][c][0.8][0]{{$\widetilde{S}_{3}$}}
\psfrag{2:1}[c][c][0.8][0]{{$2:1$}}
\epsfig{file=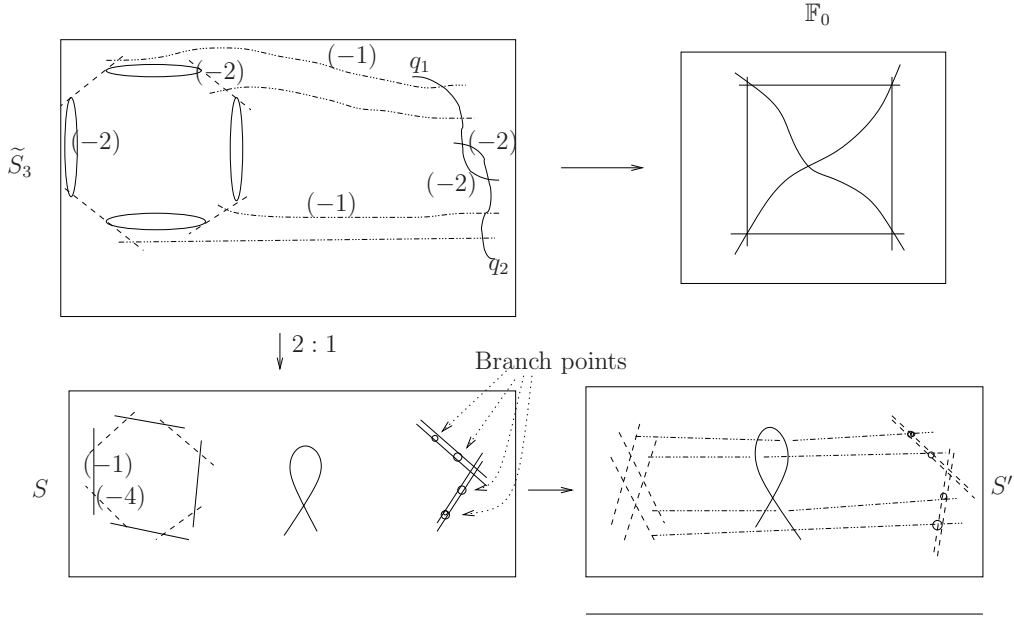,width=5.2in} 
\end{center}
\caption{Obtaining the component $\widetilde{S}_{3}$.}
\label{fig:third.component} 
\end{figure}

Finally let us look at the geometry of the third component of $Y_0$,
namely the proper transform of $Y^{\op{prop}}_0$ under the
desingularization we perform (see Figure~\ref{fig:third.component}).

\

\

As we saw $Y^{\op{prop}}_0$ is a double cover of $S'$ branched in the
$I_4$ and the double nodal fibers. There are several singular points,
besides the non-normal ones.

The blow up of $p_1\,$ and $p_2\,$ amounts to deleting the components
of the double nodal fiber from the branch locus of $Y^{\op{prop}}_0
\to S'$. But the four singular $A_1$ points on that fiber remain parts
of the branch locus.

The blow up of $n_1\cup n_2\cup n_3\cup n_4$ amounts in producing the
$I_8$-type fiber on $S'$ (with four disjoint components of
self-intersection $-4$ and four of self-intersection $-2$) and
resolving four of the singular points on $Y^{\op{prop}}_0$. Only the
$-4$-components of this fiber belong to the branch locus of the
cover. The completion of this process produces the surface $S$ with
one $I_8$-like fiber and one nodal fiber with four singular points on
it.

So, the proper preimage $\widetilde{S}_{3}$ of $Y^{\op{prop}}_0$ (the
third components of $Y_0$) is a double cover of $S$ branched in four
lines: the four $(-4)$-lines from the $I_8$-like fiber, and in four
singular points, in the nodal fiber. It is easy to realize this
surface as a blow up of $\mathbb{F}_0$ in eight points as follows.

Choose, as it was in the case of $D_1$ and $D_2$, two $(1,
1)$-divisors intersecting diagonally a quadrilateral made by four
general rulings on $\mathbb{F}_0$. Blow up the surface in the four
verices of the quadrilateral considered as schemes of depth two on the
diagonal divisors. Each divisor gets blown up four times in
infinitesimally near points to the vertices of the quadrilateral. This
is how the third component $\widetilde{S}_{3}$ of $Y_0$ looks
like. The diagonal divisors will be the curves of intersection of $F$
with $D_1$ and $D_2$ , $q_1$ and $q_2$ respectively.

It is easy to see that the points on $q_i$ where $-1$-curves intersect
them on one component match with the analogous points on the second
component to which $q_i$ belongs. This fact shows that the choices of
the rulings and diagonal divisors for each individual component of
$Y_0$ have to be compatible with each other.  \ \hfill $\Box$

\

\hfill

\noindent
From this analysis one immediately gets

\begin{cor} \label{cor:quadrics} The zero-dimentional part of the
  critical set of 
  $\bw : Y \to \mathbb{C}$ consists of two points. The rest of the
  critical set consists of three rational curves which intersect in
  two points as shown in Figure~\ref{fig-DbW} below.
\end{cor}

\begin{figure}[!ht]
\begin{center}
\psfrag{Q}[c][c][1][0]{{$q$}}
\psfrag{E1}[c][c][1][0]{{$E_{1}$}}
\psfrag{E2}[c][c][1][0]{{$E_{2}$}}
\epsfig{file=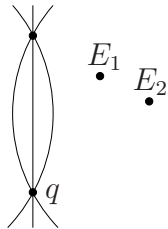,width=0.8in} 
\end{center}
\caption{The critical points of $\bw$ for the intersection of two
  quadrics in $\mathbb{P}^{5}$.}
\label{fig-DbW} 
\end{figure}

\section{The large volume comparison of Fukaya categories}
  \label{sec:large}

With Theorem~\ref{theo-genus2-symmetric} and
Theorem~\ref{theo-quadrics} at our disposal we can now use Homological
Mirror Symmetry to attack Conjecture~\ref{con:side.A.equiv} when $G =
C$ is a genus two curve and $F = X_{4}$ is a Fano threefold of degree
four.

Recall that if $C$ is a smooth complex curve of genus two, then $C$
gives rise to a pencil of quadrics in $\mathbb{P}^{5}$ whose base
locus is a smooth Fano threefold $X_{4}$ of degree four. Explicitly
since $C$ is hyperelliptic we have a double cover $C \to
\mathbb{P}^{1}$ whose branch points in the appropriate affine
coordinate on $\mathbb{P}^{1}$ are given by some complex numbers
$c_{1}$, \ldots $c_{6}$. If $(x_{1}:x_{2}:\cdots :x_{6})$ are 
homogeneous coordinates on $\mathbb{P}^{5}$, then we can associate
with $C$ the pencil of quadrics in $\mathbb{P}^{5}$ spanned by 
$\sum_{i = 1}^{6} x_{i}^{2}$ and $\sum_{i=1}^{6} c_{i}x_{i}^{2}$. The
degree four Fano threefold corresponding to $C$ is the complete
intersection 
\[
X_{4} \ : \   \sum_{i = 1}^{6} x_{i}^{2} = 0, \qquad \sum_{i=1}^{6}
c_{i}x_{i}^{2} = 0.
\]
By  \cite[Theorem~2.7]{BO} the derived category of coherent sheaves
of $C$ has a natural fully faithful embedding in the derived category
of coherent sheaves on $X_{4}$. More precisely $D^{b}(C)$ embeds in
$D^{b}(X_{4})$ as
the orthogonal of two line bundles:
\[
D^{b}(X_{4}) = \left\langle D^{b}(C), \mathcal{O}_{X_{4}},
\mathcal{O}_{X_{4}}(1) \right\rangle.
\]  
In fact \cite[Theorem~2.7]{BO} gives an explicit description of the
kernel object $K_{\Phi} \in D^{b}(C\times X_{4})$ defining the
inclusion $\Phi : D^{b}(C) \hookrightarrow D^{b}(X_{4})$. Therefore
$C$ and $X_{4}$ constitute a geometric pair that satisfies the
hypotheses of Conjecture~\ref{con:side.A.equiv}. The conjecture then
predicts that for any choice of a complexified K\"{a}hler structure
$\alpha_{F}$ on $\underline{X}_{\, 4}$ there exists an associated
choice of a complexified K\"{a}hler structure $\alpha_{G}$ on
$\underline{C}$, and a way to interpret $K_{\Phi}$ as a coisotropic
brane on $\left(\underline{C}\times \underline{X}_{\, 4},
-p_{C}^{*}\alpha_{G} + p_{X_{4}}^{*}\alpha_{F}\right)$ which induces a
fully faithful functor
\[
\Psi \ : \ \overline{\Fuk}(C,\alpha_{G}) \hookrightarrow
\overline{\Fuk}(X_{4},\alpha_{F}) 
\]
of Fukaya categories.

Before we can address the question of existence of $\Psi$ we need to
properly understand teh matching of the complexified K\"{a}hler
classes $\alpha_{F}$ and $\alpha_{G}$. The idea is to use the mirror
maps which identify certain K\"{a}hler moduli of $C$ and $X_{4}$ with
the complex moduli of the mirror Landau-Ginzburg models constructed in
theorems~\ref{theo-genus2-symmetric} and \ref{theo-quadrics}. This
presents a problem since the K\"{a}hler moduli relevant to the mirror
maps are the ones coming from equivariant K\"{a}hler classes on the
ambient toric varieties. Since $C \subset \mathbb{P}^{1}\times
\mathbb{P}^{1}$ and $X_{4} \subset \mathbb{P}^{5}$ we have two ambient
classes that show up as parameters for the Landau-Ginzburg mirror of
$C$, and only one parameter for the Landau-Ginzburg mirror of
$X_{4}$. In Section~\ref{ssec-genus2} we saw that only one of the two
ambient classes for $C \subset \mathbb{P}^{1}\times \mathbb{P}^{1}$
shows up as a modulus of the complex Landau-Ginzburg mirror of
$C$. Namely we showed that:
\begin{itemize}
\item If $B + i\omega$ is a complexified K\"{a}hler class on
$\mathbb{P}^{1}\times \mathbb{P}^{1}$, then the $B +
i\omega$-symplectic volumes $a_{1}$ and $a_{2}$ of the two rulings
appear as deformation parameters for the Landau-Ginzburg mirror of
$(\underline{C},B+i\omega)$  but only
$a_{2}$ is a modulus of the formal neighborhood of the singularities
of the zero fiber of the superpotetial.
\item If $\alpha_{F}$ is a complex multiple of the Fubini-Studi form
  on $\mathbb{P}^{5}$, then the $\alpha_{F}$-symplectic volume of the
  hyperplane in $\mathbb{P}^{5}$ appears as a deformation parameter
  for the Landau-Ginzburg mirror of $(\underline{X}_{\,
  4},\alpha_{F})$  but this parameter gives an iso-trivial deformation
  of the central fiber and so is not a modulus of the Landau-Ginzburg
  model. 
\end{itemize}
This shows that we have a one parameter ambiguity for what
$\alpha_{G}$ can be, whereas $\alpha_{F}$ is unique up to
an isomorphism. To remove this ambiguity we will consider the large
volume limit of the curve of genus $2$. Specifically we will look at
the situation where $|a_{2}| \to +\infty$. We have the following 

\begin{theo} \label{theo:main} Suppose that the Homological Mirror
  Symmetry conjecture holds for the Fano variety $X_{4}$, i.e. suppose
  that the Karoubi closure of the derived Fukaya category of
  $(\underline{X}_{\, 4}, 2\pi i \omega_{FS})$ is equivalent to the
  Karoubi closure of the derived category of the mirror
  Landau-Ginzburg model. Then Conjecture~\ref{con:side.A.equiv} holds
  for $C$ and $X_{4}$ in the large volume limit. Specifically there is a
  fully faithful functor
\[
\Psi \ : \ \lim_{\left|\int_{\mathbb{P}^{1}\times \{ \op{pt} \}} B
\right| \to -\infty} \overline{\Fuk}(\underline{C}, B+i\omega)
\hookrightarrow \overline{\Fuk}(\underline{X}_{\, 4}, 2\pi i
\omega_{FS}).
\]
\end{theo}
{\bfseries Proof:} Part of the statement of the theorem is to make
sense of limiting Fukaya category on the left hand side. We will use
mirror symmetry to show the existence of the limit and to construct
the embedding $\Psi$ in terms of the mirror Landau-Ginzburg models.

Using Seidel's proof \cite[Theorem~1.1]{SEID} of the Homological
Mirror Conjecture for $C$ we can identify the category
$\overline{\Fuk}(\underline{C},B+i\omega)$ with the Karoubi closure
$\overline{D^{b}_{\op{sing}}}(\bfY_{0})$ of the derived category of
singularities of the zero fiber of the mirror Landau-Ginzburg model.
By \cite[Theorem~2.10]{orlov-completion}
$\overline{D^{b}_{\op{sing}}}(\bfY_{0})$ depends only on the formal
neighborhood of the singularities in $\bfY_{0}$, while by our
Theorem~\ref{theo-genus2-symmetric} it follows that this formal
neighborhood depends on a single moduli parameter $\lambda$. In the
proof of Theorem~\ref{theo-genus2-symmetric} we identified $\lambda$
on one hand with the ratio of the roots of the quadratic polynomial
$x^{2} + x + a_{2}$ and on the other hand with the cross-ratio of the
four points on $\mathbb{P}^{1}$ necessary for constructing and gluing
together the three components of $\bf{Y}_{0}$. Since Seidel's theorem
identifies $\overline{\Fuk}(\underline{C},B+i\omega)$ with the Karoubi
closure of the derived category of singularities of the $\bfY_{0}$ of
modulus $\lambda$ and the limit $\lambda \to -1$ corresponds to
$|a_{2}| \to +\infty$ we see that $\lim_{|a_{2}| \to +\infty}
\overline{\Fuk}(\underline{C}, B+i\omega)$ is naturally identified
with the Karoubi closure of the derived category of singularities of
the $\bfY_{0}$ corresponding to the value $\lambda =
- 1$ of the modulus parameter. This gives a rigorous meaning of the
limiting category $\lim_{|a_{2}| \to +\infty}
\overline{\Fuk}(\underline{C}, B+i\omega)$.

Similarly in the proof of Theorem~\ref{theo-quadrics} we argued that
the zero fiber $Y_{0}$ of the Landau-Ginzburg mirror of
$(\underline{X}_{\, 4}, 2\pi i \omega_{FS})$  is built out of three
rational surfaces which are constructed and glued together out of four
points on $\mathbb{P}^{1}$ with cross-ratio $(-1)$. Furthermore in 
Corollary~\ref{cor:quadrics} and
Corollary~\ref{cor:genus.two} we saw that the critical loci of the
Landau-Ginzburg mirrors of a genus two curve and an intersection of
two quadrics in $\mathbb{P}^{5}$ are scheme-theoretically
isomorphic and that the combinatorial structures
of the full formal neighborhoods of these critical loci in the fibers
of the respective superpotentials are also the same. In particular it
folllows that the formal neighborhood of the
singularities in the $\bfY_{0}$ corresponding to the value $\lambda =
- 1$ of the modulus parameter is isomorphic to the formal neighborhood
of the singularities in the zero fiber $Y_{0}$ of the Landau-Ginzburg
mirror $(Y,\bw)$ of  $(\underline{X}_{\, 4}, 2\pi i
\omega_{FS})$. By Corollary~\ref{cor:quadrics} and
\cite[Theorem~2.10]{orlov-completion} the category
$\overline{D^{b}_{\op{sing}}}(Y_{0})$ embeds as an admissible
  subcategory in $\overline{D^{b}}(Y,\bw)$ which by our assumption is
    mirror equivalent to the Karoubi closed Fukaya category $\overline{\Fuk}
(\underline{X}_{\, 4}, 2\pi i\omega_{FS}))$. This completes the proof
    of the theorem. \ \hfill $\Box$

\

\noindent
The fact that the Fukaya category of $X_{4}$ has no moduli suggests
that the previous theorem should hold without taking the large
volume limit of the Fukaya categories of $C$. In other words, we
expect that for every $B+i\omega$ on $C$ we have a fully faithful
functor
\[
\Psi_{a_{2}} \ : \ \overline{\Fuk}(\underline{C}, B+i\omega)
\hookrightarrow \overline{\Fuk}(\underline{X}_{\, 4}, 2\pi i
\omega_{FS})
\]
given by an explicit kernel which depends only on $K_{\Phi}$ and
the $a_{2}$ period of $B + i\omega$.

This statement again has a mirror incarnation relating the complex
Landau-Ginzburg mirrors of $(\underline{C}, B+i\omega)$ and
$(\underline{X}_{\, 4}, 2\pi i \omega_{FS})$. More precisely consider 
the zero fiber $Y_{0}$ of the Landau-Ginzburg mirror of
$(\underline{X}_{\, 4}, 2\pi i \omega_{FS})$ constructed in
Theorem~\ref{theo-quadrics} and the zero fiber $\bfY_{0}^{\lambda}$ of
the Landau-Ginzburg mirror of $(\underline{C}, B+i\omega)$ constructed
in Theorem~\ref{theo-genus2-symmetric}. Here the superscript
$\lambda$ indicate that the zero fiber of the mirror of
$(\underline{C}, B+i\omega)$ corresponds to the moduli parameter
$\lambda$ which is the ratio of the roots of $x^{2} + x + a_{2} = 0$
with $a_{2}$ being the second period of $B + i\omega$. With this
notation, the existence of $\Psi_{a_{2}}$ can be recast in mirror
terms as the existence for all $\lambda$ of an equivalence
\[
\phi_{\lambda} \ : \ \overline{D^{b}_{\op{sing}}}(\bfY_{0}^{\lambda})
  \stackrel{\cong}{\to} \overline{D^{b}_{\op{sing}}}(Y_{0}).
\]

\

\noindent
In the remainder of this section we discuss a possible construction of the
functor $\Psi_{a_{2}}$ which does not use mirror symmetry and is
performed directly on the level of Fukaya categories.

For this construction we will use yet another geometric interpretation
of $X_{4}$. If $C$ is a smooth curve of genus two, the
Narasimhan-Ramanan theorem \cite{nr} identifies the corresponding Fano
$X_4$ with the moduli space of stable rank 2 bundle on $C$ of degree
one.  The underlying $C^{\infty}$ manifold of $\underline{X}_{\, 4}$
can then be identified with the moduli space of flat $SO(3)$
connections on $\underline{C}$. Now the observation is that for every
Lagrangian $s \subset C$ we can construct a Lagrangian $L_{s} \subset
X_4$. By definition $L_{s}$ is the moduli space of flat $SO(3)$
connections on $C$ that have trivial monodromy on the loop $s$.
The universal $SO(3)$ connection on $\underline{C}\times \underline{X}_{\, 4}$
allows us to transform connections on loops $s$ to connections on
$L_{s}$ and it is natural to expect that this assignment extends
to a fully faithful functor between the Fukaya categories. 

As a simple check on this expectation, we will show that the Floer
homology on both sides match in the following two very simple cases:

\begin{enumerate}
\item[1)] if two loops are disjoint then the corresponding
moduli spaces should be disjoint;
\item[2)] if two loops intersect once then the
corresponding moduli spaces should intersect once.
\end{enumerate}

\

\noindent
To be specific let $\bfa_1$, $\bfb_1$, $\bfa_2$, $\bfb_2$ be four
loops in $C$ which constitute a standard system of generators of
$\pi_1(C)$, subject to the relation $[\bfa_1,\bfb_1]\cdot
[\bfa_2,\bfb_2]=1$. Fix a symplectic form $\omega$ on $\underline{C}$
of unit volume. Then by \cite{nr} we can identify $\underline{X}_{4}$
with the moduli of pairs $(E,\nabla)$ where $E$ complex rank $2$
bundles on $\underline{C}$, and $\nabla$ is an irreducible connection
with curvature $-i\pi \omega \op{id}_{E}$. Every such connection is
characterized by its holonomy along the four generating loops
$\bfa_1$, $\bfb_1$, $\bfa_2$, $\bfb_2$. Moreover, the 1-chain
$[\bfa_1,\bfb_1]\cdot [\bfa_2,\bfb_2]$ (an 8-sided polygon) bounds a
disc in $C$, and the area of that disc is exactly the area of $C$,
i.e. is equal to $1$.  Therefore by Stokes formula and by the formula
for the curvature of the connection, the holonomy along this chain is
given by $\exp(-i\pi \op{id}_{E}) = - \op{id}_{E}$. Hence, the moduli
space $\underline{X}_{\, 4}$ can be identified with the space of 4-tuples
$(A_1,A_2,B_1,B_2)$ in  $SU(2)$ such that
$[A_1,B_1][A_2,B_2]=-\op{I}$, modulo conjugation.

Note that since the connection is not flat but only projectively flat,
the holonomy varies when we move the loop in its homotopy class;
however, if we move the loop by a Hamiltonian isotopy then the
holonomy does not change, because the area swept by the loop is
$0$. This is consistent with the fact we are interested in Lagrangians
in $C$ up to Hamiltonian isotopy only.

Now let us consider connections that are trivial both on $\bfa_1$ and
on $\bfa_2$. Since we need to have $A_1=A_2=\op{I}$, it follows that
$[A_1,B_1]\cdot [A_2,B_2]=[\op{I},B_1]\cdot [\op{I},B_2]=\op{I}$ is
never equal to $-\op{I}$. Hence there are no solutions and the two
Lagrangians $L_{\bfa_1} = \{A_1=Id\}$ and $L_{\bfa_2} = \{A_2=Id\}$ are
disjoint. On this example we have $HF(L_{\bfa_{1}},L_{\bfa_{2}}) =
HF(\bfa_1,\bfa_2) = 0$.

Similarly consider the Lagrangians in $X_{4}$ defined by say $\bfa_1$
and $\bfb_1$.  We need quadruples with $A_1=B_1=\op{I}$, so the
intersection of $L_{\bfa_{1}}$ with $L_{\bfb_1}$ consists of the
conjugacy classes of pairs $(A_2,B_2)$ in $SU(2)$ such that $
[A_2,B_2]=-\op{I}$. Up to conjugation we can diagonalize $A_2$ in the
form $diag(\exp(it),\exp(-it))$. Here $t$ is not multiple of $\pi$
since otherwise $A_2$ would be central. Then, writing the equation
$[A_2,B_2]=-Id$ explicitly in terms of coefficients we get that
\begin{itemize}
\item $B_2$  is anti-diagonal and 
\item $A_2 = \begin{pmatrix} i & 0 \\ 0 & -i \end{pmatrix}$  
or $ A_2 = \begin{pmatrix} -i & 0 \\ 0 & i \end{pmatrix}$.
\end{itemize}

\

\noindent
Since all such pairs $(A_2,B_2)$ are conjugate the intersection
of $L_{\bfa_{1}}$ with $L_{\bfb_{1}}$ consists of a single point, the same
as the intersection of $\bfa_1$ with $\bfb_1$. In other words Floer
homologies have rank 1 in both cases.

\begin{rem} It is an interesting question to understand the geometric
  nature of additional Lagrangian objects in the Fukaya category of
  $X_4$.  A natural candidate for one of them is the Lagrangian of all
  connections that come as pullbacks from the quotient of $C$ by an
  (orientation-reversing) involution.
\end{rem} 

\begin{rem} The geometric relation between $C$ and $X_4$ is not
  special to the genus two case.  Mukai \cite{MUK} showed that there
are many other Fano threefolds that can be realized as Brill-Noether
loci in the moduli space of bundles over some curve of high genus.  We
expect that an anlogue of Conjecture~\ref{con:side.A.equiv} statement
should hold in these cases as well. In paticular we expect that
Conjecture~\ref{con:side.A.equiv} should hold when $F$ is the Fano
threefold $X_{16}$ and $G$ is a curve of genus $3$, or when $F$ is a
Fano threefold $X_{12}$ and $G$ is a curve of genus 7. Some
calculations for the mirrors of reductive Fano varieties give strong
evidence for that.
\end{rem}

\begin{rem}
One can also describe the connection between Fukaya category of $S$ -
a generic intersection of quadric and cubic in $\CP^5$ and the
hypersurface $N$ of degree six in weigthed $\CP^4$ with weights
$(1,1,1,2,3)$. The affine mirror of $S$ is threefold:
\[
Y^{\op{aff}}_{S} \ : \ \left\{
\begin{split}
x_1 + x_2  & = -1 \\
x_3 + x_4 + x_5 & = -1 \\
x_1\cdot x_2 \cdot \ldots \cdot x_6 & =A 
\end{split}
\right.
\] 
equipped with the superpotential $\bw^{\op{aff}}_{S} = x_6$.

Similarly we can describe the mirror of $N$.  
Consider the potential given by $\bw^{\op{aff}}_{N} =
u_1 + \cdots + u_4 + u_5 + v$ over the fivefold $Y^{\op{aff}}_{N}$ in
${\mathbb C}^{*5} 
\times {\mathbb C}$ given by the equation:
\[
Y^{\op{aff}}_{N} \ : \ u_1 \cdot  u_2 \cdot  u_3 \cdot u_4^2 \cdot
u_5^3 = a \cdot v^6. 
\] 
Even though the fibers of $\bw^{\op{aff}}_{S}$ are two dimensional and
the fibers of $\bw^{\op{aff}}_{N}$ are four dimensional we expect 
that the critical sets of $\bw^{\op{aff}}_{S}$  and $\bw^{\op{aff}}_{N}$
have the same topology. Moreover, following an algebro geometric conjecture
made by A.Kuznetsov \cite{KUZ} and relating $D^{b}(S)$ and $D^{b}(N )$,
we expect that Conjecture~\ref{con:side.A.equiv} holds with $F = N$
and $G = S$.
\end{rem}

\section{Derived categories of singularities and flops}
\label{sec:flops} 

In this section we fill in the details of a statement that we
frequently used in our analysis: the fact that the derived category of
a potential on a three-dimensional  smooth variety does not change
under a simple flop. This statement follows from a more general
derived equivalence of Landau-Ginzburg models that we procced to
describe.

Let $Y\to S$ and $Y'\to S$ be two quasi-projective flat schemes over a
quasi-projective scheme $S$.  Let  $\E\in D^{b}(Y\times_S Y')$ be an
object such that:
$\E$ has finite Tor-amplitude over $Y$, finite
Ext-amplitude over $Y'$, and $\supp(\E)$ is projective over $Y$ and
$Y'$. Under these assumptions the integral transform  $\Phi_{\E} := R
p_{Y'*} (L p^*_{Y}(\bullet)\otimes\E)$  sends $D^{b}(Y)$
to $D^{b}(Y')$.  The right adjoint to $\Phi_{\E}$ is defined by the
formula $R p_{Y*}R\hom(\E, p^{!}_{Y'}(\bullet)),$ where $p^{!}_{Y'}$ is
the right adjoint to $p_{Y'*}.$ Since $\E$ has finite Ext-dimension over
$Y'$ the right adjoint sends $D^{b}(Y')$ to $D^{b}(Y)$.
Moreover, in this case the functor $\Phi_{\E}$ sends the subcategory
of perfect complexes $\perf(Y)$ to the subcategory of perfect complexe
$\perf(Y')$.

Let $i: s\hookrightarrow S$ be a point. Consider the fibers $i_Y:
Y_s\hookrightarrow Y$ and $i_Y':Y'_s\hookrightarrow Y'$ and denote by
$\E_s$ the object of $\db(Y_s\times Y'_s)$ which is the inverse
image $L i^*_{Y\times_S Y'} \E$ with respect to embedding
$i_{Y\times_S Y'}: Y_s\times Y'_s\hookrightarrow Y\times_S Y'.$ It can
be easily checked that the object $\E_s$ also has finite Tor-amplitude
over $Y_s,$ finite Ext-amplitude over $Y'_s$ and its support is
projective over $Y_s$ and $Y'_s$ as well.  Hence we obtain the functor
$\Phi_{\E_s}:\db(Y_s)\to \db(Y'_s)$, which sends
$\perf(Y_s)$ to $\perf(Y'_s)$ and has a right adjoint.

A direct calculation gives $R i_{Y'*}\Phi_{\E_s}\cong \Phi_{\E}R
i_{Y*}.$ Using this relation one checks immediately that if the
functor $\Phi_{\E}$ is an equivalence then the functor $\Phi_{\E_s}$
is an equivalence as well.

Note also that if the functor
$\Phi_{\E_s}:\db(Y_s)\stackrel{\sim}{\to}\db(Y'_s)$ is an equivalence,
then it induces equivalence between categories of perfect complexes,
because the perfect complexes are characterized as the homologically
finite objects in the derived category of coherent sheaves,
i.e. objects $A$ for which $R\Hom(A,B)$ has only finitely many
nontrivial cohomologies for any $B$. Thus any equivalence
$\Phi_{\E_s}$ induces an equivalence between triangulated categories
of singularities $D^{b}_{\op{sing}}(Y_s)$ and $D^{b}_{\op{sing}}(Y'_s)$.

We apply this observation to our situation. Let $\bw: Y\to \bC$ be a
3-dimensional Landau-Ginzburg-model. And let $C\in Y_0$ be a rational
curve in the fiber over $0$ such that the normal bundle to $C$ in $Y$
is isomorphic to $\mathcal{O}(-1)\oplus \mathcal{O}(-1)$. Let us
consider the flop of $Y$ in $C$. More precisely, in this situation we
can consider the blow up $\widetilde{Y}$ of $C\subset Y$ and then blow
down the exceptional divisor $C\times \bP^1$ to $\mathbb{P}^{1}$,
i.e. blow down along the ruling $C$. We obtain another
Landau-Ginzburg-model $\bw':Y'\to\bC.$ Let us take as $\E$ the
structure sheaf $\mathcal{O}_{\widetilde{Y}}$ as an object of
$\db(Y\times_{\bC}Y')$. It is shown in \cite{BO} that the functor
$\Phi_{\E}$ is an equivalence between derived categories of coherent
sheaves on $Y$ and $Y'$.  Since $Y$ and $Y'$ are smooth, any object on
the product has finite Tor-amplitude over $Y$, and finite
Ext-amplitude over $Y'$.  Moreover, the support of $\E$ is projective
over $Y$ and $Y'$. Thus the previous discussion above applies and we
get that the fibers $Y_0$ and $Y'_0$ have equivalent derived
categories of coherent sheaves and equivalent triangulated categories of
singularities. In other words we have just proven the following:

\begin{prop}
Let $\bw: Y\to\bC$ and $\bw': Y'\to\bC$ be two Landau-Ginzburg-models
which are related to each other via flop in a rational curve $C$ from
the fiber $Y_0$ over 0 whose normal bundle is isomorphic to
$\mathcal{O}(-1)\oplus\mathcal{O}(-1).$ Then the triangulated categories
of singularities $\dsing(Y_0)$ and $\dsing(Y'_0)$ are equivalent.
\end{prop}

\section{Computations with mirror branes} \label{sec:match}

In this section we return to the HMS question and
identify objects in the derived category of the
Landau-Ginzburg mirror of a genus two curve which are natural
candidates for mirrors of the standard basis of loops on the curve.

Let now $Y_0$ be the fiber of a Landau-Ginzburg model $\bw:Y\to \bC$ which is a
union of three rational surfaces $S_i$ where $i=1,2,3.$ We denote by
$F_i$ the complement to $S_i$ in $Y_0$ which is the union of $S_j$ and
$S_k,$ where $j\ne k\ne i.$ Since the relative canonical sheaf
$\omega_{Y/\bC}$ is trivial. We have that the restriction of the line
bundle $\mathcal{O}_Y(S_i)$ on $S_i$ gives us the canonical sheaf
$\omega_{S_i}.$ Note that the canonical sheaf $\omega_{S_i}$ also
coincides with the restriction of $\mathcal{O}_Y(-F_i).$ Let $C_i$ be the
intersection $S_i\cap F_i.$ It is a curve on $S_i,$ which consists of
two rational components meeting in two points and is a divisor in the
anticanonical system on $S_i.$

Consider the short exact sequence on $Y$
$$
0\to \mathcal{O}_Y(-S_i)\to \mathcal{O}_Y \to \mathcal{O}_{S_i}\to 0.
$$
It induces $\bZ/2$-periodic resolution of $\mathcal{O}_{S_i}$ as the
sheaf on $Y_0$ 
\begin{equation}\label{resol}
\{\cdots\to \mathcal{O}_{Y_0}\to \mathcal{O}_{Y_0}(-S_i)\to
\mathcal{O}_{Y_0}\}\to \mathcal{O}_{S_i}\to 0
\end{equation}

Let $L$ and $M$ be two line bundles on $Y_0.$ Denote by $L_i$ and
$M_i,$ $i=1,2,3,$ the coherent sheaves on $Y_0$ which are obtained by
restrictions of $L$ and $M$ on $S_i,$ respectively.  Tensoring the
resolution (\ref{resol}) with the line bundle $L$ we obtain a
$\bZ/2$-periodic resolution for the coherent sheaf $L_i$
\begin{equation}\label{resol2}
\{\cdots\to L\to L(-S_i)\to L\}\to L_i\to 0
\end{equation}

Using this resolution and the fact that the restriction of $O(S_i)$ on
$S_i$ coincides with the canonical sheaf $\omega_{S_i}$ we can
calculated $\rhom(L_i, M_j)$ on $Y_0.$ When $j=i$ we get the following
complex
$$ 0\to L_i^*\otimes M_i\stackrel{0}{\to} L_i^*\otimes M_i\otimes
\omega_{S_i}\stackrel{u}{\to} L_i^*\otimes M_i\stackrel{0}{\to}\cdots
$$ 
here $u$ is obtained from the injection $\mathcal{O}_{S_i}(-F_i)\to
\mathcal{O}_{S_i}$ by tensoring with $L_i^*\otimes M_i.$

When $j\ne i$ we get the complex
$$ 
0\to L_j^*\otimes M_j\stackrel{v}{\to} L_j^*\otimes
M_j(Q_{ij})\stackrel{0}{\to} L_j^*\otimes M_j\stackrel{v}{\to}\cdots
$$ 
where $v$ is obtained from the injection $\mathcal{O}_{S_j}\to
\mathcal{O}_{S_j}(Q_{ij})$ by tensoring with $L_j^*\otimes M_j.$

Summarizing we get

\begin{prop} Let $L$ and $M$ be two line bundles on $Y_0.$ And let $L_i$ and
$M_i,$ $i=1,2,3,$ be the coherent sheaves on $Y_0$ which are obtained
by restrictions of $L$ and $M$ on $S_i,$ respectively. Then
\begin{itemize}
\item[1)] when $j=i$ we have $\ext^{2p+1}_{Y_0}(L_i, M_i)\cong0$ and
$\ext^{2p+2}_{Y_0}(L_i, M_i)\cong (L^*\otimes M)|_{C_i},$

\item[2)] when $j\ne i$ we have $\ext^{2p}_{Y_0}(L_i, M_j)\cong0$ and
$\ext^{2p+1}_{Y_0}(L_i, M_j)\cong (L_j^*\otimes
M_j)(Q_{ij})|_{Q_{ij}}.$
\end{itemize}
\end{prop}

Consider now the triangulated category of singularities $\dsing{Y_0}.$
By Proposition 1.21 of \cite{O} we know that Ext's between objects in
$\dsing{Y_0}$ can be calculated as Ext's in $\coh(Y_0).$ More
precisely, in our case we have that
$$ \Hom_{\dsing{Y_0}}(L_i, M_j[N])\cong \op{Ext}^{N}_{Y_0}(L_i, M_j)
$$ for $N>2.$ But since shift by $[2]$ in $\dsing{Y_0}$ is identity,
we can calculated all Hom's in triangulated categories of
singularities of $Y_0$ taking Ext's in category of coherent sheaves on
it. Considering spectral sequence from local to global Ext's, which is
degenerated in our case we immediately obtain the following

\begin{prop}\label{ext}
Let $L$ and $M$ be two line bundles on $Y_0.$ And let $L_i$ and $M_i,$
$i=1,2,3,$ be the coherent sheaves on $Y_0$ which are obtained by
restrictions of $L$ and $M$ on $S_i,$ respectively. Then
\begin{itemize}
\item[1)] when $j=i$ we have
$$
\begin{array}{l}
\Hom_{\dsing{Y_0}}(L_i, M_i)\cong H^0(C_i, L^*\otimes M)|_{C_i})\\
\Hom_{\dsing{Y_0}}(L_i, M_i[1])\cong H^1(C_i, L^*\otimes M)|_{C_i}),
\end{array}
$$
\item[2)] when $j\ne i$ we have
$$
\begin{array}{l}
\Hom_{\dsing{Y_0}}(L_i, M_j)\cong H^1(Q_{ij}, (L_j^*\otimes
M_j)(Q_{ij})|_{Q_{ij}})\\ \Hom_{\dsing{Y_0}}(L_i, M_j[1])\cong
H^0(Q_{ij}, (L_j^*\otimes M_j)(Q_{ij})|_{Q_{ij}}).
\end{array}
$$
\end{itemize}
\end{prop}

\begin{cor}
For any $L_i$ as above we have
$$ 
\Hom_{\dsing{Y_0}}(L_i, L_i)=\bC\quad\text{and}\quad \Hom(L_i,
L_i[1])=\bC
$$
\end{cor}
This follows from the fact that $H^0(C_i, \mathcal{O}_{C_i})=\bC$ and
$H^1(C_i, \mathcal{O}_{C_i})=\bC.$

Let us take $L_i$ and $M_j,$ when $i\ne j.$ Consider the restriction
of these sheaves on $Q_{ij}\cong \bP^1.$ Let $L_i|_{Q_{ij}}\cong
O(l_i)$ and $M_j|_{Q_{ij}}\cong O(m_j).$ Then Proposition \ref{ext}
and the property that the normal bundle to $Q_{ij}$ in $S_j$ is
isomorphic to $\mathcal{O}_{\bP^1}(-1)$ imply that
$$
\Hom_{\dsing{Y_0}}(L_i, M_j)=H^1(\bP^1, O(m_j-l_i-1))\quad\text{and}\quad
\Hom_{\dsing{Y_0}}(L_i, M_j[1])=H^0(\bP^1, O(m_j-l_i-1)).
$$

On each surface $S_i$ there are two (-1)-curves each of which meet the
curve $Q_{ij}$ in one point. Denote them as $E_{ij}'$ and $E_{ij}''.$
Consider a set of four sheaves $\mathcal{O}_{S_1}(E_{12}'),
\mathcal{O}_{S_1}(E_{13}'), \mathcal{O}_{S_2}, \mathcal{O}_{S_3}$ on
$Y_0.$ Then nontrivial Hom's in the triangulated category of
singularities between them are only the following
$$
\Hom_{\dsing{Y_0}}(\mathcal{O}_{S_1}(E_{12}'),
\mathcal{O}_{S_2})=\bC,\quad \Hom_{\dsing{Y_0}}(\mathcal{O}_{S_2},
\mathcal{O}_{S_1}(E_{12})'[1])=\bC, 
$$
$$
\Hom_{\dsing{Y_0}}(\mathcal{O}_{S_1}(E_{13}'),
\mathcal{O}_{S_3})=\bC,\quad \Hom_{\dsing{Y_0}}(\mathcal{O}_{S_3},
\mathcal{O}_{S_1}(E_{13})'[1])=\bC. 
$$

\begin{cor}
Hom's between objects $\mathcal{O}_{S_1}(E_{12}'),
\mathcal{O}_{S_1}(E_{13}'), \mathcal{O}_{S_2}, \mathcal{O}_{S_3}$ in
triangulated category of singularities of $Y_0$ coincide with Floer
cohomologies between Lagrangian circles on a curve of genus two which
represent the standard basis $\bfa_1, \bfa_2, \bfb_1, \bfb_2$ in the first
integral homologies of the curve.
\end{cor}

\section{Appendix}

In this appendix we collect some important technical facts on the
structure of the singular fibers of the Landau-Ginzburg mirrors of a
curve of genus two or the intersection of two quadrics in
$\mathbb{P}^{5}$. 

\subsection{Local analytic structure at the curves of intersection}

Our description of the component $\widetilde{S}_{3}$ of the zero fiber
of either the mirror of $C$ or the mirror of $X_{4}$ shows that this
surface has an elliptic fibration in which one $I_{2}$ fiber is given
by the union of the curves $q_1 \cup q_2$. This is why the formal or
analytic neighborhoods of $q_{1}\cup q_{2}$ can be identified in the
two mirrors. The global elliptic fibrations on the $\widetilde{S}_{3}$
components are different though: in the case of the genus two curve
the fibration has two fibers of type $I_1$, and one of type $I_2^*$
while in the case of the intersection of two quadrics in $\CP^5$ there
are two fibers of type $I_1$ and one of type $I_8$.

The same applies to the curve $q_3$. In our birational models for the
mirrors they are singular (nodal) fibers in conic fibrations, and
hense have the same analytic neighbourhoods.  In order to identify the
full analytic neighborhods of the singularities of the central fibers
for the two mirrors we have to identify the respective curves $q_3$ so
that {\bf four} points on them match. Since the cross ratios of the
four tuple of points in the two mirrors are different ($\lambda$ and
$(-1)$ respectively) it is not a priori clear if this can be
done geometrically. We still expect that the corresponding categories
are equivalent though.

\subsection{Symmetries in the central fibers}

Here is one more feature the fibrations we are working with have.

As mentioned in the previous sections, the pairs diagonal $(1,
1)$-divisors in $\mathbb{F}_0$ and the pairs of positive sections in
$\mathbb{F}_2$ have to be chosen so that certain agreements are
met. This follows from the symmetry (coming basically from the initial
equations of the varieties $Y^{\op{aff}}$) according to which the
exceptional components of the central fibers are isomorphic (this is
before we blow up the last two points as explained before, and via an
involution of the central fibers which leave the third component
invariant).

In the case of $\mathbb{F}_0$ corresponding to the component
$\widetilde{S}_{3}$, this implies that the pair of diagonals have to
get exchanged by an involution of $\mathbb{F}_0$ acting so that maps
one of the members of the one pair of rulings to the other one and
leaves the members of the other pair invariant. This actually forces
to have the same property on both other components as well. This in
turn gives us that the three components form an $S_3$ symmetric
configuration up to some extra blow-ups.

The same can be said about the sections of $\mathbb{F}_2$. In a sense,
this reflects the symmetry we had on the general fiber related to the
three involutions acting there.

\newcommand{\etalchar}[1]{$^{#1}$}

\end{document}